\documentclass[11pt]{article}

\usepackage[T1]{fontenc}
\usepackage{lmodern}
\usepackage{microtype}
\usepackage{amsmath,amssymb,amsthm,mathtools}
\usepackage{bm}
\usepackage{enumitem}

\usepackage[margin=1.15in]{geometry}
\usepackage{graphicx}

\usepackage{hyperref}
\usepackage[nameinlink,capitalise,noabbrev]{cleveref}

\hypersetup{
    colorlinks=true,
    linkcolor=blue,
    citecolor=blue,
    urlcolor=blue
}

\newtheorem{theorem}{Theorem}[section]
\newtheorem{lemma}[theorem]{Lemma}

\theoremstyle{definition}
\newtheorem{definition}[theorem]{Definition}

\theoremstyle{remark}

\title{A Class of Coupled Forward-Backward Systems Under General Filtrations}
\author{
    Matthew Keyworth\\
    Department of Statistics, University of Warwick\\
    \texttt{matthew.keyworth@warwick.ac.uk}
}
\date{}

\begin{document}

\maketitle

\begin{abstract}
We consider a class of coupled-forward backward systems defined on a general probability space assumed only to satisfy the usual conditions. We obtain existence and uniqueness results for arbitrary finite time horizons.
\end{abstract}

%\tableofcontents
\section{Introduction}
In this paper we examine a class of coupled forward-backward systems defined on a general probability space assumed only to satisfy the usual conditions. Such equations arise for example in stochastic control from the stochastic maximum principle (see for example \cite{YongBook:1999} for the case of coupled forward-backward stochastic differential equations). We study questions of existence and uniqueness of solutions without assuming any connection to an underlying stochastic control problem. These questions have been studied extensively in the literature in the case of coupled FBSDEs (see for example \cite{MaBook:2007}, \cite{ZhangBook:2017}). For general coupled forward-backward systems existence and uniqueness can typically be obtained for sufficiently small time horizons via contraction arguments (see for example \cite{AntonelliPaper:1993}). For arbitrary time horizons existence and uniqueness in the case of coupled FBSDEs with Brownian driver can be obtained assuming certain monotonicity and global Lipschitz conditions on the coefficients (see \cite{PengPaper:1995}, \cite{YongPaper:1997}, \cite{PengPaper:1999}, \cite{PardouxPaper:1999}). Alternatively, 
the so called four-step scheme enables one to solve certain coupled FBSDEs by representing the solution of the backward equation in terms of a solution to an associated PDE (see \cite{ProtterPaper:1994}). In our case established FBSDE theory does not apply directly since, although we can ensure that the backward equation admits a c\'adl\'ag version with a semimartingale decomposition, we cannot assume that the martingale part admits a stochastic integral representation. Consequently we take a different approach, and by lifting the problem to a suitable Hilbert space we prove that solutions of our coupled forward-backward system are related to the fixed points of a certain nonlinear operator. As such, studying existence and uniqueness of solutions equates to studying the fixed points of this nonlinear operator. Uniqueness of solutions for all finite time horizons can be established under only certain monotonicity assumptions on the coefficients. Proving existence for arbitrary time-horizons proves to be more challenging but can be established by first proving existence for a suitably chosen affine system and then proving an appropriate continuation lemma in a similar way to the so called method of continuation, which can be used to establish solvability in certain cases for coupled FBSDEs (see \cite{YongPaper:1997} and \cite{MaBook:2007}). In \cref{sec:preliminaries} we recall some basic facts and establish our notation. In \cref{sec:problem_setup} we define our coupled forward-backward system of interest and state our technical assumptions. In \cref{sec:main_results} we detail the main results of the paper and how they fit together. In \cref{sec:proofs_main_res} we work through proofs of our main results.

\section{Preliminaries and Notation} \label{sec:preliminaries}
We first recall some basic facts and introduce the notation used throughout. Define the following sets of stochastic processes:
\begin{equation}
    \begin{aligned}
         \mathcal{L}([0,T]; \mathbb{R}^d) &:=  \left\{\text{$X:\Omega\times[0,T]\mapsto \mathbb{R}^d: \ X$  is $\mathcal{F}_T\otimes\mathcal{B}([0,T])$-measurable}\right\} \\
         \mathcal{L}_{\mathbb{F}}([0,T]; \mathbb{R}^d) &:= \left\{\text{$X \in \mathcal{L}([0,T]; \mathbb{R}^d)$ and $\mathbb{F}$-adapted}\right\} \\
        \mathcal{L}_{\mathbb{F}}^D([0,T]; \mathbb{R}^d) &:= \left\{X \in  \mathcal{L}_{\mathbb{F}}([0,T]; \mathbb{R}^d) : \text{$t\mapsto X_t$ c\`adl\`ag $\mathbb{P}\text{-a.s.}$}\right\}\\
        \mathcal{L}_{\mathbb{F}}^C([0,T]; \mathbb{R}^d) &:= \left\{X \in  \mathcal{L}_{\mathbb{F}}^D([0,T]; \mathbb{R}^d) : \text{$t\mapsto X_t$ continuous $\mathbb{P}\text{-a.s.}$}\right\} \\
        \mathcal{M}^2([0,T]; \mathbb{R}^d) &:= \left\{X \in \mathcal{L}^{D}_{\mathbb{F}}([0,T]; \mathbb{R}^d): \text{ $X$ is an $\mathbb{F}$-martingale with } \mathbb{E}[||X_T||_2^2] < \infty\right\} 
    \end{aligned}
\end{equation}
For any subset $A\subseteq \mathcal{L}([0,T]; \mathbb{R}^d)$ and any pseudo-metric $d:A\times A\mapsto \mathbb{R}_+$ using the equivalence relation:
\[X\sim_d Y \iff d(X, Y) = 0, \quad X,Y \in A\]
we may define the equivalence classes:
\begin{align*}
    [X]_d := \{Y \in A: Y\sim_d X\}, \quad \forall X \in A
\end{align*}
The pseudo-metrics:
\begin{align*}
    d_p(X, Y) := \mathbb{E}\left[\int_0^T ||X_t- Y_t ||_p^p\,ds\right]^{1/p}, \quad d_S(X, Y) = \mathbb{E}\left[\sup_{t \in [0,T]}||X_t - Y_t||_2^2\right]^{1/2}
\end{align*}
induce $d\mathbb{P}\otimes dt\text{-a.e.}$ equivalence and equivalence up to indistinguishability respectively. We thus define the following sets of equivalence classes:
\begin{equation}\label{eq:process_spaces}
    \begin{aligned}
        \mathbb{L}^p_{\mathbb{F}}([0,T]; \mathbb{R}^d) &:= \left\{[X]_{d_p} : X \in \mathcal{L}_{\mathbb{F}}([0,T]; \mathbb{R}^d), d_p(X, 0) <\infty\right\}, \quad p \in \mathbb{N} \\
        \mathbb{S}^2_{\mathbb{F}}([0,T]; \mathbb{R}^d) &:= \left\{[X]_{d_S}: X \in \mathcal{L}^D_{\mathbb{F}}([0,T]; \mathbb{R}^d), d_S(X, 0) < \infty\right\} \\
        \mathbb{S}^{2,c}_{\mathbb{F}}([0,T]; \mathbb{R}^d) &:= \left\{[X]_{d_S}: X \in \mathcal{L}^C_{\mathbb{F}}([0,T]; \mathbb{R}^d), d_S(X, 0) < \infty\right\} \\
         \mathbb{M}^2_{\mathbb{F}}([0,T]; \mathbb{R}^d) &:= \left\{[X]_{d_S}: X \in \mathcal{M}^2([0,T]; \mathbb{R}^d), \ d_S(X, 0) < \infty \right\} 
    \end{aligned}
\end{equation}
For a set of equivalence classes $\mathbb{A}$ with respect to the equivalence relation $\sim_d$, we will use the usual abuse of notation and take $X \in \mathbb{A}$ to mean a representer of the equivalence class $[X]_d \in \mathbb{A}$.
We recall some basic facts about the spaces that will be used in the subsequent sections. Equipped with the inner product:
\begin{align*}
    \langle X, Y \rangle &:= \mathbb{E}\left[\int_0^T\langle X_t, Y_t\rangle\,dt\right]
\end{align*}
the space $\mathbb{L}^2_{\mathbb{F}}([0,T]; \mathbb{R}^d)$ is a Hilbert space. Equipped with the norm induced by $d_S$, the spaces $\mathbb{S}^{2,c}_{\mathbb{F}}([0,T]; \mathbb{R}^d)$ and $\mathbb{S}^{2}_{\mathbb{F}}([0,T]; \mathbb{R}^d)$ are Banach spaces. We make use of canonical continuous embedding $\mathbb{S}^{2}_{\mathbb{F}}\xhookrightarrow{i}\mathbb{L}^2_{\mathbb{F}}$ where:
\begin{align*}
    i:\mathbb{S}^2_{\mathbb{F}}([0,T]; \mathbb{R}^d)\rightarrow \mathbb{L}^2_{\mathbb{F}}([0,T]; \mathbb{R}^d), \quad i([X]_{d_S}) = [X]_{d_2}
\end{align*}
and denote by $\mathcal{A}$ its image $\mathcal{A} := i(\mathbb{S}^{2}_{\mathbb{F}})\subseteq \mathbb{L}^2_{\mathbb{F}}$ then we have:
\begin{align} \label{eq:def_mathcal_A}
    \mathcal{A} &= \{[X]_{d_2} \in \mathbb{L}^2_{\mathbb{F}}([0,T]; \mathbb{R}^d): \text{$[X]_{d_2}$ admits an $\mathbb{F}$-adapted representative that is c\`adl\`ag $\mathbb{P}\text{-a.s.}$}\}
\end{align}
Since $\mathcal{A}$ is a set of equivalence classes, from now on we implicitly fix a representative. For $X \in \mathcal{A}$ we select the unique (up to indistinguishability) $\mathbb{F}$-adapted $\mathbb{P}\text{-a.s.}$ c\`adl\`ag representative. By virtue of the usual conditions each $\mathbb{F}$ adapted martingale admits a c\`adl\`ag modification. Consequently for any square integrable martingale $M$ we can take $M \in \mathbb{M}^2_{\mathbb{F}}([0,T]; \mathbb{R})$.
\section{Problem Formulation and Assumptions} \label{sec:problem_setup}
With $T \in (0,\infty)$ let $(\Omega, \mathcal{F}, \mathbb{F}:=(\mathcal{F}_t)_{t \in [0,T]}, \mathbb{P})$ be a filtered probability space satisfying the usual conditions. Let $g:\Omega \times \mathbb{R}^d\rightarrow \mathbb{R}^d$, $b:\Omega\times [0,T] \times \mathbb{R}^d\rightarrow \mathbb{R}^d$ and $\bar{\rho}:\Omega\times [0,T]\times \mathbb{R}^d\rightarrow U$ for convex $U\subseteq \mathbb{R}^d$. As is the usual convention of notation we mostly suppress the dependence of functions $g, b, \rho$ on $\omega$. We study existence and uniqueness of solutions to the following coupled forward-backward system:
\begin{equation} \label{eq:fbs_general}
    \begin{aligned}
        v_t &= v_0 + \int_0^t\bar{\rho}(s,\lambda_s)\,ds, \quad t \in [0,T], \quad v_0 \in \mathbb{R}^d\\
        \lambda_t &= \mathbb{E}\left[g(v_T) + \int_t^Tb(s, v_s)\,ds \ \bigg| \mathcal{F}_t \right], \quad \mathbb{P}\text{-a.s.}, \quad t \in [0,T] 
    \end{aligned}
\end{equation}
\begin{definition}[Solution Concept]\label{def:solution_concept}
We say that any pair $(v, \lambda) \in \mathbb{S}_{\mathbb{F}}^{2,c}([0,T]; \mathbb{R}^d)\times \mathbb{S}^2_{\mathbb{F}}([0,T]; \mathbb{R}^d)$ satisfying \cref{eq:fbs_general} is a solution.  We say that the forward-backward system has a unique solution iff for any two solutions $(v,\lambda)$ and $(v',\lambda')$ we have that:
\begin{align*}
    \mathbb{P}(v_t = v'_t; 0\leq t\leq T) = 1, \quad \mathbb{P}(\lambda_t = \lambda'_t; 0\leq t\leq T)= 1
\end{align*}
\end{definition}
We examine questions of existence and uniqueness of solutions to \cref{eq:fbs_general} under the following assumptions:
\paragraph{Integrability Assumptions.}
With $v^\rho_t = v_0 + \int_0^t \rho_s\,ds, \ t \in [0,T]$:
\begin{enumerate}[label=(A\arabic*),  ref=A\arabic*]
 \item \label{ass:A1} $g(v_T^\rho) \in L^2(\Omega, \mathcal{F}_T, \mathbb{P}; \mathbb{R}^d), \ \forall \rho \in \mathbb{L}^2_{\mathbb{F}}([0,T]; \mathbb{R}^d)$
    \item \label{ass:A2} $b^\rho \in \mathbb{L}^2_{\mathbb{F}}([0,T]; \mathbb{R}^d), \forall \rho \in \mathbb{L}^2_{\mathbb{F}}([0,T]; \mathbb{R}^d)$ with $b^\rho_t := b(t,v^\rho_t), \ t \in [0,T]$
\end{enumerate}
\paragraph{Monotonicity Assumptions.}
\begin{enumerate}[label=(M\arabic*),  ref=M\arabic*]
    \item \label{ass:B1} $\langle g(\omega, v) - g(\omega, v'), v - v'\rangle \geq 0, \quad \forall v, v' \in \mathbb{R}^d, \quad \mathbb{P}\text{-a.s.}$
    \item \label{ass:B2} $\exists \xi \in \mathbb{L}^1_{\mathbb{F}}([0,T]; \mathbb{R}), \ \xi_t > 0, \ \mathbb{P}\text{-a.s.}, \text{ a.e. } t \in [0,T]$ such that: 
    \[\langle b(t,\omega, v)-b(t,\omega, v'), v-v'\rangle \geq \xi(t,\omega)||v - v'||_2^2, \quad \forall v, v' \in \mathbb{R}^d, \ \mathbb{P}\text{-a.s.}, \text{ a.e. } t \in [0,T]\]
    \item \label{ass:B3}  $\exists\mu > 0$ such that $\forall \lambda, \lambda' \in \mathbb{R}^d:$
    \begin{align*}
        &\langle\bar{\rho}(t,\omega, \lambda) - \bar{\rho}(t,\omega, \lambda'), \lambda - \lambda'\rangle  \leq -\mu||\bar{\rho}(t,\omega, \lambda) - \bar{\rho}(t,\omega, \lambda')||_2^2, \quad \mathbb{P}\text{-a.s.}, \ \text{a.e.} \ t \in [0,T]  
    \end{align*}
\end{enumerate}
\paragraph{Global Lipschitz Assumptions.}
\begin{enumerate}[label=(A\arabic*), start=3, ref=A\arabic*]
    \item \label{ass:A5} There exist constants $L_g, L_b > 0$ such that:
    \begin{align*}
        ||g(\omega, v) - g(\omega, v')||_2 &\leq L_g||v - v'||_2, \quad v,v' \in \mathbb{R}^d, \quad \mathbb{P}\text{-a.s.} \\
        ||b(t,\omega, v) - b(t,\omega, v')||_2 &\leq L_b||v - v'||_2, \quad v,v' \in \mathbb{R}^d, \quad \mathbb{P}\text{-a.s.}, \ \text{a.e. } t \in [0,T]
    \end{align*}
\end{enumerate}
Note that by \ref{ass:B3} with $L_\rho = 1/\mu$ we obtain the global Lipschitz condition:
\begin{align*}
    ||\rho(t,\omega, \lambda) - \rho(t,\omega, \lambda')||_2 &\leq L_\rho||\lambda - \lambda'||_2, \quad \lambda, \lambda' \in \mathbb{R}^d, \quad \mathbb{P}\text{-a.s.}, \ \text{a.e. } t \in [0,T] 
\end{align*}

\section{Main Results} \label{sec:main_results}
Here we relay the main results of the paper and how they fit together. 
\paragraph{Decoupled Systems and a Duality Relation}
Defining the the set:
\begin{align*}
    \mathcal{U}[0,T]:=\{\rho \in \mathcal{L}_{\mathbb{F}}([0,T]; \mathbb{R}^d): \ \rho(t,\omega) \in U, \ \mathbb{P}\text{-a.s.} \ \text{ a.e } t \in [0,T]\}
\end{align*}
which modulo $d\mathbb{P}\otimes dt\text{-a.e.}$ equivalence becomes the convex subset $\mathbb{U}[0,T]\subseteq \mathbb{L}^2_{\mathbb{F}}([0,T]; \mathbb{R}^d)$, we embed the coupled forward backward system of \cref{eq:fbs_general} into a large set of decoupled forward-backward systems related to a strictly monotone operator via a duality pairing. In particular, examine the set of decoupled forward-backward systems:
\begin{equation} \label{eq:weakly_coupled_fbs}
    \begin{aligned}
        v_t &= v_0 + \int_0^t\rho_s\,ds, \quad \rho \in \mathbb{U}[0,T] \\
        \lambda_t^\rho &= \mathbb{E}\left[g(v_T) + \int_t^Tb(s,v_s)\,ds \ \bigg|\mathcal{F}_t\right], \quad \mathbb{P}\text{-a.s.}, \quad t \in [0,T]    
    \end{aligned}
\end{equation}
and define the strictly monotone operator $A:\mathbb{L}^2_{\mathbb{F}}([0,T]; \mathbb{R}^d)\rightarrow \mathbb{L}^2_{\mathbb{F}}([0,T]; \mathbb{R}^d)^*$:
\begin{equation}\label{eq:strictly_monotone_operator}
    \begin{aligned}
        A(\rho)(\varphi) &:= \mathbb{E}\left[g(v_T)\cdot v_T' + \int_0^Tb(t,v_t)\cdot v_t'\,dt\right] \\
        v_t &= v_0 + \int_0^t\rho_s\,ds, \quad  t \in [0,T], \quad \rho \in \mathbb{L}^2_{\mathbb{F}}([0,T]; \mathbb{R}^d)\\
        v_t' &= \int_0^t\varphi_s\,ds, \quad t \in [0,T], \quad \varphi \in \mathbb{L}^2_{\mathbb{F}}([0,T]; \mathbb{R}^d)
    \end{aligned}
\end{equation}
We first prove some useful intermediate results. Namely, that $\lambda^\rho \in \mathcal{A}, \ \forall \rho \in \mathbb{U}[0,T]$, and that the following duality relation holds:
\begin{equation} \label{eq:duality_rel}
    \begin{aligned}
        A(\rho)(\varphi) = \langle \lambda^\rho, \varphi\rangle_{\mathbb{L}^2_{\mathbb{F}}}, \quad \forall \varphi \in \mathbb{L}^2_{\mathbb{F}}([0,T]; \mathbb{R}^d), \ \forall \rho \in \mathbb{U}[0,T]
    \end{aligned}
\end{equation}
I.e. for each $\rho \in \mathbb{U}[0,T]$, $\lambda^\rho \in \mathcal{A}$ is the Riesz representative of the continuous linear functional $A(\rho) \in \mathbb{L}^2_{\mathbb{F}}([0,T]; \mathbb{R}^d)^*$. 

\paragraph{The Coupled System as a Fixed Point Problem}
Establishing these results crucially enables us to re-frame the coupled forward-backward system of \cref{eq:fbs_general} as a fixed point problem on the Hilbert space $\mathbb{L}^2_{\mathbb{F}}([0,T]; \mathbb{R}^d)$. We can define a Nemytskii map $\mathcal{N}:\mathbb{L}^2_{\mathbb{F}}([0,T]; \mathbb{R}^d)\rightarrow \mathbb{U}[0,T]$ such that:
\begin{equation}\label{eq:nemytskii_map}
    \begin{aligned}
        \mathcal{N}(\lambda)(t,\omega) = \bar{\rho}(t,\omega, \lambda(t,\omega)), \quad (t,\omega) \in [0,T]\times \Omega, \quad \lambda \in \mathbb{L}^2_{\mathbb{F}}([0,T]; \mathbb{R}^d)    
    \end{aligned}
\end{equation}
and from \ref{eq:weakly_coupled_fbs} write down the associated set of decoupled forward-backward systems as follows:
\begin{equation} \label{eq:decoupled_fbs_parameterised}
    \begin{aligned}
        v_t &= v_0 + \int_0^t\mathcal{N}(\lambda)_s\,ds, \quad t \in [0,T], \quad \lambda \in \mathbb{L}^2_{\mathbb{F}}([0,T]; \mathbb{R}^d) \\
        \lambda_t^{\mathcal{N}(\lambda)} &= \mathbb{E}\left[g(v_T)+ \int_t^Tb(s,v_s)\,ds \ \bigg| \mathcal{F}_t\right], \quad \mathbb{P}\text{-a.s.}, \quad t \in [0,T] 
    \end{aligned}
\end{equation}
Using the duality relation of \ref{eq:duality_rel} we would have that:
\begin{align}\label{eq:duality_rel_coupled}
    A(\mathcal{N}(\lambda))(\varphi) = \langle  \lambda^{\mathcal{N}(\lambda)}, \varphi \rangle_{\mathbb{L}^2_{\mathbb{F}}}, \quad \forall \varphi \in \mathbb{L}^2_{\mathbb{F}}([0,T]; \mathbb{R}^d), \ \lambda \in \mathbb{L}^2_{\mathbb{F}}([0,T]; \mathbb{R}^d)
\end{align}
where $\lambda^{\mathcal{N}(\lambda)} \in \mathcal{A}$. Now denoting the Riesz map by $R:\mathbb{L}^2_{\mathbb{F}}([0,T]; \mathbb{R}^d)\mapsto \mathbb{L}^2_{\mathbb{F}}([0,T]; \mathbb{R}^d)^*$, we can define an important self map $\Phi:\mathcal{A}\rightarrow \mathcal{A}$ such that:
\begin{align} \label{eq:non_lin_phi_map}
    \lambda^{\mathcal{N}(\lambda)} &= \Phi(\lambda) := (R^{-1}\circ A\circ \mathcal{N})(\lambda), \quad \lambda \in \mathcal{A}
\end{align}
Defining the set of fixed points $\operatorname{Fix}\Phi := \{\lambda \in \mathcal{A} : \lambda = \Phi(\lambda)\}$, we will prove in \cref{theorem:sol_iff_fixed_point} that the pair $(v^\lambda,\lambda) \in \mathbb{S}^{2,c}_{\mathbb{F}}([0,T]; \mathbb{R}^d)\times \mathbb{S}^{2}_{\mathbb{F}}([0,T]; \mathbb{R}^d)$ solves \cref{eq:fbs_general} iff $\lambda \in \operatorname{Fix}\Phi$. Therefore finding a solution to the coupled forward-backward system of \cref{eq:fbs_general} can be done by solving the fixed point problem to obtain $\lambda^* \in \operatorname{Fix}\Phi$ and then feeding $\lambda^*$ into the forward equation to obtain the solution pair $(v^*, \lambda^*)$. The conceptual re-framing from coupled forward-backward system to fixed point problem eventually leads to our main results of existence and uniqueness of solutions.
\paragraph{Monotonicity and Uniqueness of Solutions}
In \cref{theorem:fbs_unique_sol} we prove that if $\operatorname{Fix}\Phi$ is non-empty, then it is a singleton, leading to our first main result: uniqueness of solutions to \cref{eq:fbs_general}, where these results follow directly from our monotonicity assumptions.

\paragraph{Existence of Solutions}
Existence of solutions to \cref{eq:fbs_general} proves to be more challenging. First we prove existence and uniqueness of solutions for an appropriately chosen affine system in \cref{subsec:existence_affine_system}. We then prove a suitable continuation lemma in \cref{lemma:continuation_lemma}, relying on a crucial a priori estimate proved in \cref{lemma:key_apriori_estimate}. These results combined eventually lead us to our second main result proved in \cref{theorem:fbs_general_existence}: existence of solutions to \cref{eq:fbs_general}.

\section{Proofs of Main Results} \label{sec:proofs_main_res}
\subsection{A Fixed Point Problem and Uniqueness of Solutions}\label{subsec:fp_and_uniqueness}
In this section we assume that: \ref{ass:A1}, \ref{ass:A2}, \ref{ass:B1}, \ref{ass:B2} and \ref{ass:B3} hold. 
\begin{lemma}\label{lemma:lambda_cadlag_sm}
     Recall the decoupled forward-backward system of \cref{eq:weakly_coupled_fbs}. Assuming \ref{ass:A1} and \ref{ass:A2}, for each $\rho \in \mathbb{U}[0,T]$ we have that 
     \begin{enumerate}
         \item $\lambda^\rho \in \mathcal{A}$ and has the semimartingale decomposition
           \begin{align*}
            \lambda_t^\rho = -\int_0^tb(s,v_s)\,ds + M_t^\rho, \quad \mathbb{P}\text{-a.s.},\ t \in [0,T], \quad \lambda_T = g(v_T), \quad \mathbb{P}\text{-a.s.}
        \end{align*}
        where $M^\rho \in \mathbb{M}^{2}_{\mathbb{F}}([0,T]; \mathbb{R}^d)$.
        \item For every $\rho \in \mathbb{U}[0,T]$ we have:
         \begin{align*}
            ||v^\rho_t||_2 &\leq ||v_0||_2 + \sqrt{T}||\rho||_{\mathbb{L}^2_{\mathbb{F}}}, \ \forall t \in [0,T] \\
            ||v^\rho||_{\mathbb{L}^2_{\mathbb{F}}} &\leq \sqrt{T}||v_0||_2 + T||\rho||_{\mathbb{L}^2_{\mathbb{F}}} 
         \end{align*}
     \end{enumerate}
    \begin{proof}
        Taking $\rho \in \mathbb{U}[0,T]$, then for each $t \in [0,T]$ by Cauchy-Schwarz and Jensen's inequality:
        \begin{align*}
            \mathbb{E}[||v_t||_2^2]^{1/2} = \mathbb{E}[||v_0 + \int_0^t\rho_s\,ds||_2^2]^{1/2} \leq ||v_0||_2 + \sqrt{T}||\rho||_{\mathbb{L}^2_{\mathbb{F}}} < \infty, \quad \forall t \in [0,T]
        \end{align*}
        thus we have $v_t \in L^2(\Omega, \mathcal{F}_t, \mathbb{P}; \mathbb{R}^d)$. We have that $v \in \mathbb{L}^2_{\mathbb{F}}([0,T]; \mathbb{R}^d)$ :
        \begin{align*}
            ||v||_{\mathbb{L}^2_{\mathbb{F}}}  = \mathbb{E}\left[\int_0^T||v_0 + \int_0^t\rho_s\,ds||_2^2\,dt\right]^{1/2} \leq \sqrt{T}||v_0||_2 + T||\rho||_{\mathbb{L}^2_{\mathbb{F}}} < \infty
        \end{align*}
        Now defining:
         \begin{align*}
            M_t^\rho &= \mathbb{E}\left[g(v_T) + \int_0^Tb(s,v_s)\,ds \ \bigg| \mathcal{F}_t\right], \quad \mathbb{P}\text{-a.s.}, \ t \in [0,T] \\
            v_t &= v_0 + \int_0^t\rho_s\,ds, \quad t \in [0,T]
        \end{align*}
        We have that by \ref{ass:A1} and \ref{ass:A2}:
        \begin{align*}
            \mathbb{E}[||M_T^\rho||_2^2]^{1/2} &= \mathbb{E}[||g(v_T) + \int_0^Tb(s,v_s)\,ds||_2^2]^{1/2} \\
            &\leq \mathbb{E}[||g(v_T)||_2^2]^{1/2} + \sqrt{T}\mathbb{E}\left[\int_0^T||b(s,v_s)||_2^2\,ds\right]^{1/2} < \infty
        \end{align*}
        Thus $M^\rho$ is a square integrable martingale with a c\`adl\`ag version by the corollary of Theorem 3.1 in \cite{LipsterBook1:2001}. From now on taking $M^\rho \in \mathbb{M}^{2}_{\mathbb{F}}([0,T]; \mathbb{R}^d)$ we have that:
        \begin{align*}
            \lambda_t^\rho = -\int_0^tb(s,v_s)\,ds + M_t^\rho, \quad \mathbb{P}\text{-a.s.},\ t \in [0,T], \quad \lambda_T = g(v_T), \quad \mathbb{P}\text{-a.s.}
        \end{align*}
        thus $\lambda^\rho$ admits a c\`adl\`ag version. Taking such a c\`adl\`ag version we have that: 
        \begin{align*}
            ||\lambda_t^\rho||_2^2 \leq 2T\int_0^T||b(s,v_s)||_2^2\,ds + 2\sup_{0\leq t\leq T}||M^\rho_t||_2^2, \quad t \in [0,T]
        \end{align*}
        Taking expectations and by Doob's inequality we obtain:
        \begin{align*}
            ||\lambda^\rho||_{\mathbb{L}^2_{\mathbb{F}}}^2 = \mathbb{E}\left[\int_0^T||\lambda^\rho_t||_2^2\,dt\right] \leq 2T^2\mathbb{E}\left[\int_0^T||b(s,v_s)||_2^2\,ds\right] + 2T\mathbb{E}\left[\sup_{0\leq t\leq T}||M^\rho_t||_2^2\right] \\
            \leq 2T^2\mathbb{E}\left[\int_0^T||b(s,v_s)||_2^2\,ds\right] + 8T\mathbb{E}\left[||M^\rho_T||_2^2\right] < \infty
        \end{align*}
        and hence $\lambda^\rho \in \mathcal{A}$
    \end{proof}
\end{lemma}

\begin{lemma}\label{lemma:A_duality_pairing} 
    Recall the decoupled forward-backward system of \cref{eq:weakly_coupled_fbs} and the operator $A:\mathbb{L}^2_{\mathbb{F}}([0,T]; \mathbb{R}^d)\rightarrow\mathbb{L}^2_{\mathbb{F}}([0,T]; \mathbb{R}^d)^*$ as defined in \ref{eq:strictly_monotone_operator}. The duality relation of \cref{eq:duality_rel} holds: for each $\rho \in \mathbb{U}[0,T]$ $\lambda^\rho \in \mathcal{A}$ is the Riesz representative of the continuous linear functional $A(\rho) \in \mathbb{L}^2_{\mathbb{F}}([0,T]; \mathbb{R}^d)^*$.
    \begin{proof}
       Taking $\rho \in \mathbb{U}[0,T]$ and taking $\mathcal{S}^{b,c\`ag}_{\mathbb{F}}([0,T]; \mathbb{R}^d)$ to be the set of simple $\mathbb{F}$-adapted c\`agl\`ad essentially bounded processes (which are dense in $\mathbb{L}^2_{\mathbb{F}}([0,T]; \mathbb{R}^d)$ by \cref{lemma:simple_proc_dense_LpF}):
        \begin{align*}
            \mathcal{S}^{b,c\`ag}_{\mathbb{F}}([0,T]; \mathbb{R}^d):=\{\varphi:[0,T]\times \Omega \rightarrow \mathbb{R}^d: \varphi(t,\omega) := \xi_0(\omega)1_{\{0\}}(t) + \sum_{i=1}^{n-1}\xi_i(\omega)1_{(t_i, t_{i+1}]}(t), \\  \xi_i \in L^\infty(\Omega, \mathcal{F}_{t_i}, \mathbb{P}; \mathbb{R}^d), \ 0 = t_1 < t_2 < \dots, t_n = T, \ n \in \mathbb{N}\}
        \end{align*}
         then for any $\varphi \in \mathcal{S}_{\mathbb{F}}^{b,c\`ag}$, we have that $v'_{\cdot} := \int_0^{\cdot}\varphi_s\,ds$ is bounded continuous. By \cref{lemma:lambda_cadlag_sm} $\lambda^\rho \in \mathcal{A}$ and admits a semimartingale decomposition. Thus using integration by parts for semimartingales we obtain:
        \begin{align*}
            g(v_T)\cdot v_T' = \int_0^T\lambda_s^\rho\cdot\varphi_s\,ds -\int_0^Tb(t, v_t)\cdot v'_t \,dt + \int_0^T v_t'\cdot dM_t^\rho
        \end{align*}
         with the stochastic integral a true martingale. Now taking expectations and rearranging, we obtain the duality pairing:
        \begin{align*}
            A(\rho)(\varphi) = \mathbb{E}\left[g(v_T)\cdot v_T' + \int_0^T b(t, v_t)\cdot v_t' \,dt \right] = \mathbb{E}\left[\int_0^T\lambda_t^\rho\cdot \varphi_t\,dt\right] = \langle \lambda^\rho, \varphi\rangle_{\mathbb{L}^2_{\mathbb{F}}}, \quad \forall \varphi \in \mathcal{S}_b
        \end{align*}  
        Where the RHS is defined for $\lambda^\rho$ up to $d\mathbb{P}\otimes dt\text{-a.s.}$ equivalence. Now by Cauchy-Schwarz we have that:
        \begin{align*}
            |A(\rho)(\varphi)| = |\langle \lambda^\rho, \varphi\rangle_{\mathbb{L}^2_{\mathbb{F}}}| \leq ||\lambda^\rho||_{\mathbb{L}^2_{\mathbb{F}}}||\varphi||_{\mathbb{L}^2_{\mathbb{F}}} \implies  \sup_{\{\varphi : ||\varphi||_{\mathcal{S}_b}\leq 1\}}|A(\rho)(\varphi)| \leq ||\lambda^\rho||_{\mathbb{L}^2_{\mathbb{F}}} < \infty
        \end{align*}
        and thus $A(\rho)$ is a continuous linear functional on $\mathcal{S}^{b,c\`ag}_{\mathbb{F}}([0,T]; \mathbb{R}^d)$ by lemma 5.1 in \cite{AltBook:2016}. Now we can extend to a continuous linear functional on $\mathbb{L}^2_{\mathbb{F}}([0,T]; \mathbb{R}^d)$. For any $\varphi \in \mathbb{L}^2_{\mathbb{F}}([0,T]; \mathbb{R}^d)$ take an approximating sequence $(\varphi_n)_{n \in \mathbb{N}}\subseteq \mathcal{S}_{\mathbb{F}}^b$ such that $ ||\varphi_n - \varphi||_{\mathbb{L}^2_{\mathbb{F}}}\rightarrow 0$, then we have the pairing:
    \begin{align*}
        A(\rho)(\varphi_n) = \langle \lambda^\rho, \varphi_n\rangle_{\mathbb{L}^2_{\mathbb{F}}}, \quad \forall n \in \mathbb{N}
    \end{align*}
    and we have by lemma 5.1 in \cite{AltBook:2016} that there exists $C > 0$ such that:
    \begin{align*}
        |A(\rho)(\varphi_n) - A(\rho)(\varphi)| = |A(\rho)(\varphi_n - \varphi)| \leq C||\varphi_n - \varphi||_{\mathbb{L}^2_{\mathbb{F}}} \rightarrow 0,\ \text{ as } n\rightarrow \infty
    \end{align*}
    We also have by Cauchy-Schwarz:
    \begin{align*}
        |\langle \lambda^\rho, \varphi_n\rangle_{\mathbb{L}^2_{\mathbb{F}}} - \langle \lambda^\rho, \varphi\rangle_{\mathbb{L}^2_{\mathbb{F}}}| = |\langle \lambda^\rho, \varphi_n - \varphi\rangle_{\mathbb{L}^2_{\mathbb{F}}} | \leq ||\lambda^\rho||_{\mathbb{L}^2_{\mathbb{F}}}||\varphi_n - \varphi||_{\mathbb{L}^2_{\mathbb{F}}} \rightarrow 0,\ \text{ as } n\rightarrow \infty
    \end{align*}
    Hence:
    \begin{align*}
        A(\rho)(\varphi) = \lim_{n\rightarrow\infty}A(\rho)(\varphi_n) = \lim_{n\rightarrow\infty}\langle \lambda^\rho, \varphi_n\rangle_{\mathbb{L}^2_{\mathbb{F}}} = \langle \lambda^\rho, \varphi\rangle_{\mathbb{L}^2_{\mathbb{F}}}
    \end{align*}
    and the duality pairing holds on $\mathbb{L}^2_{\mathbb{F}}([0,T]; \mathbb{R}^d)$:
    \begin{align*}
         A(\rho)(\varphi) = \langle \lambda^\rho, \varphi\rangle_{\mathbb{L}^2_{\mathbb{F}}}, \quad \forall \varphi \in \mathbb{L}^2_{\mathbb{F}}([0,T]; \mathbb{R}^d), \quad \forall \rho \in \mathbb{U}[0,T]
    \end{align*}
    By the Riesz representation theorem this duality pairing is unique and $\lambda^\rho \in \mathcal{A}$ is the Riesz representative of $A(\rho)$. 
    \end{proof}
\end{lemma}
\begin{lemma}\label{lemma:strictly_monotone_operator}
    The operator $A:\mathbb{L}^2_{\mathbb{F}}([0,T]; \mathbb{R}^d)\rightarrow\mathbb{L}^2_{\mathbb{F}}([0,T]; \mathbb{R}^d)^*$ as defined in \ref{eq:strictly_monotone_operator} is strictly monotone (and hence injective):
    \begin{proof}
        For $\rho, \rho' \in \mathbb{L}^2_{\mathbb{F}}([0,T]; \mathbb{R}^d)$  with $\rho\neq \rho'$ we have:
        \begin{align*}
            (A(\rho)-A(\rho'))(\rho - \rho')&= \mathbb{E}\left[(g(v_T) - g(v'_T))\cdot (v_T - v'_T) + \int_0^T(b(t,v_t) - b(t,v'_t))\cdot (v_t - v_t')\,dt\right] \\
            v_t &= v_0 + \int_0^t\rho_s\,ds, \quad v'_t = v_0 + \int_0^t\rho_s'\,ds, \quad t \in [0,T]
        \end{align*}
        By \ref{ass:B1} and \ref{ass:B2} we have that:
        \begin{align*}
            A(\rho)(\rho - \rho') - A(\rho')(\rho - \rho') \geq \mathbb{E}\left[\int_0^T\xi_t||v_t - v'_t||_2^2\,dt\right] \geq 0
        \end{align*}
        Clearly, since $\xi_t > 0$ $\mathbb{P}\text{-a.s.}$ $\text{a.e.}$ $t \in [0.T]$, we have that the RHS is zero if and only if $v_t - v_t'=0$ $\mathbb{P}\text{-a.s.}$ $\text{a.e.}$ $t \in [0.T]$ and by the fact that:
        \begin{align*}
            v_t - v_t' = \int_0^t(\rho_s - \rho'_s)\,ds, \quad t \in [0,T]
        \end{align*}
        we can conclude that:
        \begin{align*}
            (A(\rho)-A(\rho'))(\rho - \rho') > 0, \quad \forall \rho\neq \rho'
        \end{align*}
        and thus $A$ is strictly monotone. To prove injectivity suppose that for some $\rho'\neq \rho \in \mathbb{L}^2_{\mathbb{F}}([0,T]; \mathbb{R}^d)$ we have that $A(\rho) = A(\rho') \in \mathbb{L}^2_{\mathbb{F}}([0,T]; \mathbb{R}^d)^*$. Thus we would have:
        \begin{align*}
            (A(\rho)-A(\rho'))(\varphi) = 0, \quad \forall \varphi \in \mathbb{L}^2_{\mathbb{F}}
        \end{align*}
        But setting $\varphi = \rho - \rho'$ contradicts strict monotonicity and thus $A$ is injective.
    \end{proof}
\end{lemma}
\begin{theorem} \label{theorem:sol_iff_fixed_point}
    A pair $(v, \lambda) \in \mathbb{S}^{2,c}_{\mathbb{F}}([0,T]; \mathbb{R}^d)\times \mathbb{S}^2_{\mathbb{F}}([0,T]; \mathbb{R}^d)$ is a solution to \cref{eq:fbs_general} if and only if $\lambda \in \operatorname{Fix}\Phi\subseteq \mathcal{A}$.
    \begin{proof}
        Taking $\lambda \in \operatorname{Fix}\Phi\subseteq \mathcal{A}$ and recalling the definition of $\Phi:\mathcal{A}\rightarrow\mathcal{A}$ (see \ref{eq:non_lin_phi_map}) we have that:
        \begin{align*}
            \lambda = \Phi(\lambda) := R^{-1}(A(\mathcal{N}(\lambda)))
        \end{align*}
        equivalently we have that:
        \begin{align*}
            A(\mathcal{N}(\lambda))(\varphi) = \langle \lambda, \varphi\rangle_{\mathbb{L}^2_{\mathbb{F}}}, \quad \forall \varphi \in \mathbb{L}^2_{\mathbb{F}}([0,T]; \mathbb{R}^d)
        \end{align*}
        But by \cref{lemma:A_duality_pairing} we must have that:
        \begin{align*}
           A(\mathcal{N}(\lambda))(\varphi) = \langle \lambda^{{\mathcal{N}}(\lambda)}, \varphi\rangle_{\mathbb{L}^2_{\mathbb{F}}}, \quad \forall \varphi \in \mathbb{L}^2_{\mathbb{F}}([0,T]; \mathbb{R}^d)
        \end{align*}
        with $\lambda^{{\mathcal{N}}(\lambda)} \in \mathcal{A}$. Thus we have that
        \begin{align*}
            \langle \lambda - \lambda^{{\mathcal{N}}(\lambda)}, \varphi\rangle_{\mathbb{L}^2_{\mathbb{F}}} = 0, \quad \forall \varphi \in \mathbb{L}^2_{\mathbb{F}}([0,T]; \mathbb{R}^d)
        \end{align*}
        and hence $\lambda = \lambda^{\mathcal{N}(\lambda)} \ d\mathbb{P}\otimes dt\text{-a.e.}$ and thus are of the same equivalence class in $\mathcal{A}$. Recall that $\mathcal{A}$ was constructed as the image of the injective mapping $i$ with $\mathbb{S}^{2}_{\mathbb{F}}([0,T]; \mathbb{R}^d)\xhookrightarrow{i}\mathbb{L}^2_{\mathbb{F}}([0,T]; \mathbb{R}^d)$ such that $\mathcal{A} := i(\mathbb{S}^{2}_{\mathbb{F}}([0,T]; \mathbb{R}^d))\subseteq \mathbb{L}^2_{\mathbb{F}}([0,T]; \mathbb{R}^d)$. Thus after fixing c\`adl\`ag representatives (see \ref{eq:def_mathcal_A}) we must have:
        \begin{align*}
        \mathbb{P}\left(\lambda_t = \lambda^{\mathcal{N}(\lambda)}; \forall t \in [0,T]\right) = 1
        \end{align*}
        now recalling \cref{eq:decoupled_fbs_parameterised}, and by indistinguishability of  $\lambda$ and $\lambda^{\mathcal{N}(\lambda)}$ we must have that $(v^\lambda,\lambda) \in \mathbb{S}^{2,c}_{\mathbb{F}}([0,T]; \mathbb{R}^d)\times \mathbb{S}^2_{\mathbb{F}}([0,T]; \mathbb{R}^d)$ solves \cref{eq:fbs_general}. Consequently:
        \begin{align*}
            \lambda \in \operatorname{Fix}\Phi \implies (v^\lambda, \lambda) \in \mathbb{S}^{2,c}_{\mathbb{F}}([0,T]; \mathbb{R}^d)\times \mathbb{S}^2_{\mathbb{F}}([0,T]; \mathbb{R}^d) \text{ solves \cref{eq:fbs_general} }
        \end{align*}
        Now supposing $(v^\lambda, \lambda) \in \mathbb{S}^{2,c}_{\mathbb{F}}([0,T]; \mathbb{R}^d)\times \mathbb{S}^2_{\mathbb{F}}([0,T]; \mathbb{R}^d)$ solves \cref{eq:fbs_general}. Recalling \cref{eq:decoupled_fbs_parameterised} it must be the case that $\mathbb{P}(\lambda_t = \lambda^{\mathcal{N}(\lambda)}_t) = 1, \ \forall t \in [0,T]$. Since $\lambda, \lambda^{\mathcal{N}(\lambda)} \in \mathcal{A}$ and two c\'adl\'ag versions are indistinguishable, $\lambda, \lambda^{\mathcal{N}(\lambda)}$ are of the same equivalence class in $\mathcal{A}$ (see \ref{eq:def_mathcal_A}). Thus by \cref{lemma:A_duality_pairing} we have that:
        \begin{align*}
            A(\mathcal{N}(\lambda))(\varphi) = \langle \lambda^{\mathcal{N}(\lambda)}, \varphi\rangle_{\mathbb{L}^2_{\mathbb{F}}} = \langle \lambda, \varphi\rangle_{\mathbb{L}^2_{\mathbb{F}}}, \quad \forall \varphi \in \mathbb{L}^2_{\mathbb{F}}([0,T]; \mathbb{R}^d)
        \end{align*}
        Equivalently we have that:
        \begin{align*}
            \lambda = \lambda^{\mathcal{N}(\lambda)} = \Phi(\lambda) := R^{-1}(A(\mathcal{N}(\lambda)))
        \end{align*}
        and $\lambda \in \operatorname{Fix}\Phi$. Consequently: 
        \begin{align*}
            (v^\lambda, \lambda) \in \mathbb{S}^{2,c}_{\mathbb{F}}([0,T]; \mathbb{R}^d)\times \mathbb{S}^2_{\mathbb{F}}([0,T]; \mathbb{R}^d) \text{ solves \cref{eq:fbs_general}} \implies \lambda \in \operatorname{Fix}\Phi
        \end{align*}
    \end{proof}
\end{theorem}
\begin{theorem}[Uniqueness of Solutions] \label{theorem:fbs_unique_sol}
If the set $\operatorname{Fix}\Phi$ is non-empty, it is a singleton. Consequently any solution to \cref{eq:fbs_general} is unique up to indistinguishability.
\begin{proof}
    For $i \in \{1, 2\}$ assume that $\lambda^i \in \operatorname{Fix}\Phi$ and define $\rho^i:=\mathcal{N}(\lambda^i) \in \mathbb{U}[0,T]$.
    \begin{align*}
        v_t^i &= v_0 + \int_0^t\rho^i_s\,ds, \quad t \in [0,T]
    \end{align*}
    then by \cref{theorem:sol_iff_fixed_point} $(v^i, \lambda^i) \in \mathbb{S}^{2,c}_{\mathbb{F}}([0,T]; \mathbb{R}^d)\times \mathbb{S}^2_{\mathbb{F}}([0,T]; \mathbb{R}^d)$ solves \cref{eq:fbs_general}. Now defining $\Delta \rho := \rho^1 - \rho^2 \in \mathbb{L}^2_{\mathbb{F}}([0,T]; \mathbb{R}^d)$, by \cref{lemma:A_duality_pairing} we have that:
    \begin{align*}
        A(\rho^1)(\Delta\rho) - A(\rho^2)(\Delta\rho) - \langle \lambda^1 - \lambda^2, \Delta\rho\rangle_{\mathbb{L}^2_{\mathbb{F}}} = 0
    \end{align*}
    by \cref{lemma:strictly_monotone_operator} we have that:
    \begin{align*}
        A(\rho^1)(\Delta\rho) - A(\rho^2)(\Delta\rho) \geq \mathbb{E}\left[\int_0^T\xi_t||v_t^1 - v_t^2||_2^2\,dt\right] > 0, \quad \rho^1 \neq \rho^2
    \end{align*}
    and by \ref{ass:B3} we have:
    \begin{align*}
        \langle \lambda^1 - \lambda^2, \Delta\rho\rangle_{\mathbb{L}^2_{\mathbb{F}}} = \mathbb{E}\left[\int_0^T(\lambda_t^1 - \lambda_t^2)\cdot(\bar{\rho}(t,\lambda^1_t) - \bar{\rho}(t,\lambda^2_t))\,dt\right] \leq 0
    \end{align*}
    Thus we obtain:
    \begin{align*}
        0 &= A(\rho^1)(\Delta\rho) - A(\rho^2)(\Delta\rho) - \langle \lambda^1 - \lambda^2, \Delta\rho\rangle_{\mathbb{L}^2_{\mathbb{F}}} \\
        &\geq  A(\rho^1)(\Delta\rho) - A(\rho^2)(\Delta\rho) \\ &\geq \mathbb{E}\left[\int_0^T\xi_t||v_t^1 - v_t^2||_2^2\,dt\right] > 0, \quad \rho^1 \neq \rho^2
    \end{align*}
    and it must be the case that $\rho^1 = \rho^2$. From the inequality:
    \begin{align*}
        \mathbb{E}\left[\sup_{0\leq t\leq T}||v_t^1 - v_t^2||_2^2\,dt\right] \leq T\mathbb{E}\left[\int_0^T||\rho^1_t - \rho^2_t||_2^2\,dt\right] = 0
    \end{align*}
    we have that $\mathbb{P}(v^1_t = v^2_t; \forall t \in [0,T])= 1$ consequently it must be the case that $\mathbb{P}(\lambda^1_t = \lambda^2_t)= 1, \ \forall t \in [0,T]$. Since two c\'adl\'ag versions are indistinguishable the result follows from \cref{def:solution_concept} and \cref{theorem:sol_iff_fixed_point}.
\end{proof}
\end{theorem}

\subsection{Existence and Uniqueness for a Suitably Chosen Affine System} \label{subsec:existence_affine_system}
Let us consider the following affine coupled forward-backward system:
\begin{equation} \label{eq:fbs_affine}
    \begin{aligned}
        v_t &= v_0 + \int_0^t(\tilde{F}_s\lambda_s + \tilde{f}_s)\,ds, \quad t \in [0,T], \quad v_0 \in \mathbb{R}^d\\
        \lambda_t &= \mathbb{E}\left[(\tilde{G}v_T + \tilde{g}) + \int_t^T(\tilde{B}_sv_s + \tilde{b}_s)\,ds \ \bigg| \mathcal{F}_t \right], \quad \mathbb{P}\text{-a.s.}, \quad t \in [0,T] 
    \end{aligned}
\end{equation}
Where we assume that:
\begin{enumerate}
    \item $\tilde{G} \in \mathbb{R}^{d\times d}$ is positive semidefinite, $\tilde{B} \in \mathbb{R}^{d\times d}$ is strongly positive definite and $\tilde{F} \in \mathbb{R}^{d\times d}$ is strongly negative definite. Note that by applying \cref{lemma:firm_non_expans} to $-\tilde{F}$ we obtain \ref{ass:B3}.
    \item $\tilde{g} \in L^2(\Omega, \mathcal{F}_T, \mathbb{P}; \mathbb{R}^d)$, and $\tilde{b}, \tilde{f}\in \mathbb{L}^2_{\mathbb{F}}([0,T]; \mathbb{R}^d)$
\end{enumerate}   
Consequently we can readily apply the results of \cref{subsec:fp_and_uniqueness}. Note that since $\tilde{f} \in \mathbb{L}^2_{\mathbb{F}}([0,T]; \mathbb{R}^d)$, in this case we have that $\mathbb{U}[0,T] := \mathbb{L}^2_{\mathbb{F}}([0,T]; \mathbb{R}^d)$.
Now defining the affine operator $A:\mathbb{L}^2_{\mathbb{F}}([0,T]; \mathbb{R}^d)\rightarrow \mathbb{L}^2_{\mathbb{F}}([0,T]; \mathbb{R}^d)^*$:
\begin{equation} \label{eq:affine_A}
    \begin{aligned}
          A(\rho)(\varphi) &= \mathbb{E}\left[(\tilde{G}v_T + \tilde{g})\cdot v'_T + \int_0^T(\tilde{B}_sv_s + \tilde{b}_s)\cdot v_s'\,ds\right] \\
    v_t &= v_0 + \int_0^t\rho_s\,ds, \quad t \in [0,T], \quad \quad \rho \in \mathbb{L}^2_{\mathbb{F}}([0,T]; \mathbb{R}^d) \\
    v'_t &= \int_0^t\varphi_s\,ds, \quad t \in [0,T], \quad  \varphi \in \mathbb{L}^2_{\mathbb{F}}([0,T]; \mathbb{R}^d)
    \end{aligned}
\end{equation}
and the affine Nemytskii map $\mathcal{N}:\mathbb{L}^2_{\mathbb{F}}([0,T]; \mathbb{R}^d)\rightarrow \mathbb{L}^2_{\mathbb{F}}([0,T]; \mathbb{R}^d)$:
\begin{align} \label{eq:nemytskii_affine}
    \mathcal{N}(\lambda)(t,\omega) = \tilde{F}\lambda(t,\omega) + \tilde{f}(t,\omega), \quad (t,\omega) \in [0,T]\times \Omega
\end{align}
By applying \cref{lemma:A_duality_pairing} we obtain the duality relation:
\begin{align*}
        A(\mathcal{N}(\lambda))(\varphi) = \langle \lambda^{\mathcal{N}(\lambda)}, \varphi\rangle_{\mathbb{L}^2_{\mathbb{F}}}, \quad \forall \lambda \in \mathcal{A}
\end{align*}
and by applying \cref{theorem:sol_iff_fixed_point} we have that a pair $(v, \lambda) \in \mathbb{S}^{2,c}_{\mathbb{F}}([0,T]; \mathbb{R}^d)\times \mathbb{S}^{2}_{\mathbb{F}}([0,T]; \mathbb{R}^d)$ solves \cref{eq:fbs_affine} iff $\lambda \in \operatorname{Fix}\Phi$ where $\Phi:=R^{-1}\circ A\circ \mathcal{N}:\mathcal{A}\rightarrow\mathcal{A}$.
\begin{lemma}[Uniqueness of Solutions to the Affine System] \label{lemma:uniqueness_affine_case}
    If the set $\operatorname{Fix}\Phi$ is non-empty, it is a singleton. Consequently any solution to \cref{eq:fbs_affine} is unique up to indistinguishability.
    \begin{proof}
        The result follows directly from \cref{theorem:fbs_unique_sol}.
    \end{proof}
\end{lemma}
\begin{lemma}[Existence for the Affine System] \label{lemma:existence_affine_case}
    The set $\operatorname{Fix}\Phi\neq \emptyset$ and thus a solution to the affine coupled forward-backward system of \cref{eq:fbs_affine} exists.
    \begin{proof}  
    Define the linear operator $A_L:\mathbb{L}^2_{\mathbb{F}}([0,T];\mathbb{R}^d)\rightarrow \mathbb{L}^2_{\mathbb{F}}([0,T];\mathbb{R}^d)^*$ such that:
    \begin{align*}
        A(\rho)_L(\varphi) &= \mathbb{E}\left[\tilde{G}(v_T-v_0)\cdot v'_T + \int_0^T\tilde{B}_s(v_s-v_0)\cdot v_s'\,ds\right] \\
        v_t &= v_0 + \int_0^t\rho_s\,ds, \quad t \in [0,T], \quad \rho \in \mathbb{L}^2_{\mathbb{F}}([0,T]; \mathbb{R}^d) \\
        v'_t &= \int_0^t\varphi_s\,ds, \quad t \in [0,T], \quad  \varphi \in \mathbb{L}^2_{\mathbb{F}}([0,T]; \mathbb{R}^d)
    \end{align*}
    and define $A_c \in \mathbb{L}^2_{\mathbb{F}}([0,T];\mathbb{R}^d)^*$ such that:
    \begin{align*}
        A_c(\varphi) &= \mathbb{E}\left[(\tilde{G}v_0 + \tilde{g})\cdot v'_T + \int_0^T(\tilde{B}_sv_0 + \tilde{b}_s)\cdot v_s'\,ds\right] \\
        v'_t &= \int_0^t\varphi_s\,ds, \quad t \in [0,T],  \quad  \varphi \in \mathbb{L}^2_{\mathbb{F}}([0,T]; \mathbb{R}^d)
    \end{align*}
    then $A= A_L + A_c$ with $A$ as defined in \ref{eq:affine_A}. Defining $\mathcal{N}_L:\mathbb{L}^2_{\mathbb{F}}([0,T]; \mathbb{R}^d)\rightarrow \mathbb{L}^2_{\mathbb{F}}([0,T]; \mathbb{R}^d)$ such that:
    \begin{align*}
        \mathcal{N}_L(\lambda)(t,\omega) = \tilde{F}\lambda(t,\omega), \quad (t,\omega) \in [0,T]\times \Omega
    \end{align*}
    We have that $\mathcal{N} = \mathcal{N}_L + \tilde{f}$ with $\mathcal{N}$ as defined in \ref{eq:nemytskii_affine}. Consequently we have that:
    \begin{align*}
        \Phi(\lambda) &= R^{-1}(A(\mathcal{N}(\lambda))) \\
        &= R^{-1}(A_L(\mathcal{N}(\lambda))) + R^{-1}(A_c) \\ 
        &= R^{-1}(A_L(\mathcal{N}_L(\lambda)) + R^{-1}(A_L(\tilde{f})+A_c) \\
        &= \Phi_L(\lambda) +\Phi_c
    \end{align*}
    with:
    \begin{align*}
        \Phi_L(\lambda) :=  R^{-1}(A_L(\mathcal{N}_L(\lambda)), \quad \Phi_c:=R^{-1}(A_L(\tilde{f})+A_c)
    \end{align*}
    we have that $\Phi_L(\lambda):\mathcal{A}\rightarrow \mathcal{A}$ is linear as a composition of linear maps, and $\Phi_c \in \mathcal{A}$.
    Thus we have that for $\lambda \in \mathcal{A}$:
    \begin{align*}
        \lambda \in \operatorname{Fix}\Phi \iff (I - \Phi)(\lambda) = 0 \iff (I - \Phi_L)(\lambda) = \Phi_c \iff \Phi_c \in \operatorname{Ran}(I - \Phi_L)
    \end{align*}
    Note that $\Phi_L$ is also well defined as a mapping $\Phi_L :\mathbb{L}^2_{\mathbb{F}}([0,T]; \mathbb{R}^d)\rightarrow\mathcal{A}\subseteq\mathbb{L}^2_{\mathbb{F}}([0,T]; \mathbb{R}^d)$. We will first show that the equation $(I-\Phi_L)(\lambda) = h$ has a solution $\lambda \in \mathbb{L}^2_{\mathbb{F}}([0,T]; \mathbb{R}^d)$ for every $h \in \mathbb{L}^2_{\mathbb{F}}([0,T]; \mathbb{R}^d)$. Applying \cref{lemma:strictly_monotone_operator} we have that that $A_L$ is strictly monotone, and by \cref{lemma:A_duality_pairing} we obtain the duality relation:
    \begin{align*}
        A_L(\mathcal{N}_L(\lambda))(\varphi) = \langle \Phi_L(\lambda), \varphi\rangle_{\mathbb{L}^2_{\mathbb{F}}}, \quad \forall \varphi \in \mathbb{L}^2_{\mathbb{F}}([0,T]; \mathbb{R}^d), \ \lambda \in \mathbb{L}^2_{\mathbb{F}}([0,T]; \mathbb{R}^d)
    \end{align*}
    By monotonicity of $A_L$ and the duality relation we have that $\forall \lambda \in \mathbb{L}^2_{\mathbb{F}}([0,T]; \mathbb{R}^d)$:
    \begin{align*}
        (A_L(\mathcal{N}_L(0)) - A_L(\mathcal{N}_L(\lambda)))(\mathcal{N}_L(0) - \mathcal{N}_L(\lambda)) = \langle \Phi_L(0) - \Phi_L(\lambda), \mathcal{N}_L(0)-\mathcal{N}_L(\lambda)\rangle_{\mathbb{L}^2_{\mathbb{F}}} \geq 0
    \end{align*}
    Thus we have that:
    \begin{align*}
        \langle (I- \Phi_L)(\lambda), \mathcal{N}_L(0)-\mathcal{N}_L(\lambda)\rangle_{\mathbb{L}^2_{\mathbb{F}}}  \geq -\langle \Phi_L(0)  - \lambda, \mathcal{N}_L(0)-\mathcal{N}_L(\lambda)\rangle_{\mathbb{L}^2_{\mathbb{F}}}
    \end{align*}
    Since $\Phi_L, \mathcal{N}_L$ are linear we have that $\Phi_L(0) = 0$ and $\mathcal{N}_L(0) = 0$ thus:
    \begin{align*}
        \langle (I-\Phi_L)(\lambda), -\mathcal{N}_L(\lambda)\rangle_{\mathbb{L}^2_{\mathbb{F}}}  \geq -\langle  \lambda, \mathcal{N}_L(\lambda)\rangle_{\mathbb{L}^2_{\mathbb{F}}}
    \end{align*}
    By assumption $\tilde{F}$ is strongly negative definite, so there exists $c > 0$ such that:
    \begin{align*}
        \langle (I-\Phi_L)(\lambda), -\mathcal{N}_L(\lambda)\rangle_{\mathbb{L}^2_{\mathbb{F}}}  \geq -\langle  \lambda, \mathcal{N}_L(\lambda)\rangle_{\mathbb{L}^2_{\mathbb{F}}} \geq c^2|| \lambda||_{\mathbb{L}^2_{\mathbb{F}}}^2
    \end{align*}
    and by Cauchy-Schwarz we obtain:
    \begin{align*}
        ||(I- \Phi_L)(\lambda)||_{\mathbb{L}^2_{\mathbb{F}}}||\mathcal{N}_L(\lambda)||_{\mathbb{L}^2_{\mathbb{F}}}  \geq  c^2|| \lambda||_{\mathbb{L}^2_{\mathbb{F}}}^2, \quad \forall \lambda \in \mathbb{L}^2_{\mathbb{F}}([0,T]; \mathbb{R}^d)
    \end{align*}
    and since $\mathcal{N}_L$ is strongly continuous it is bounded by lemma 5.1 in \cite{AltBook:2016}, and we have:
    \begin{align*}
        ||(I- \Phi_L)(\lambda)||_{\mathbb{L}^2_{\mathbb{F}}}||\mathcal{N}_L||_{op}||\lambda||_{\mathbb{L}^2_{\mathbb{F}}}  \geq  c^2|| \lambda||_{\mathbb{L}^2_{\mathbb{F}}}^2, \quad \forall \lambda \in \mathbb{L}^2_{\mathbb{F}}([0,T]; \mathbb{R}^d)
    \end{align*}
    thus we obtain the inequality:
    \begin{align*}
        ||(I- \Phi_L)(\lambda)||_{\mathbb{L}^2_{\mathbb{F}}}  \geq  \frac{c^2}{||\mathcal{N}_L||_{op}}|| \lambda||_{\mathbb{L}^2_{\mathbb{F}}}, \quad \forall \lambda \in \mathbb{L}^2_{\mathbb{F}}([0,T]; \mathbb{R}^d)
    \end{align*}
    Thus, we have $(I-\Phi_L):\mathbb{L}^2_{\mathbb{F}}([0,T]; \mathbb{R}^d)\rightarrow \mathbb{L}^2_{\mathbb{F}}([0,T]; \mathbb{R}^d)$ to be injective since:
    \begin{align*}
         (I-\Phi_L)(\lambda) = 0 \implies 0 \geq \frac{c^2}{||\mathcal{N}_L||_{op}}||\lambda||_{\mathbb{L}^2_{\mathbb{F}}} \iff \lambda = 0 \implies \operatorname{Ker}(I - \Phi_L) = \{0\}
    \end{align*}
    and consequently we obtain the apriori estimate:
    \begin{align*}
        \frac{c^2}{||\mathcal{N}_L||_{op}}\operatorname{dist}(\lambda, \operatorname{Ker}(I-\Phi_L)) \leq ||(I- \Phi_L)(\lambda)||_{\mathbb{L}^2_{\mathbb{F}}}, \quad \forall \lambda \in \mathbb{L}^2_{\mathbb{F}}([0,T]; \mathbb{R}^d)
    \end{align*}
    Hence by Banach's closed range theorem (Theorem 3.E in \cite{ZeidlerBook2:1995}) $\operatorname{Ran}(I- \Phi_L)$ is closed and by Corollary 1 of the same theorem $\operatorname{Ran}(I-\Phi_L) = \operatorname{Ker}((I-\Phi_L)^*)^\perp$ where $(I-\Phi_L)^*$ denotes the Hilbert adjoint operator (see 12.2 in \cite{AltBook:2016}). Thus we have that:
    \begin{align*}
        (I-\Phi_L) \text{ surjective } \iff \mathbb{L}^2_{\mathbb{F}}([0,T]; \mathbb{R}^d)\subseteq  \operatorname{Ker}((I-\Phi_L)^*)^\perp \iff \operatorname{Ker}((I-\Phi_L)^*) = \{0\}
    \end{align*}
    Now supposing $h \in \operatorname{Ker}((I-\Phi_L)^*)$ and $h\neq 0$. We have that:
    \begin{align*}
        \langle \lambda, (I-\Phi_L)^*(h)\rangle_{\mathbb{L}^2_{\mathbb{F}}} = \langle (I- \Phi_L)(\lambda), h\rangle_{\mathbb{L}^2_{\mathbb{F}}} = 0, \quad \forall \lambda \in \mathbb{L}^2_{\mathbb{F}}([0,T]; \mathbb{R}^d)
    \end{align*}
    Since by assumption $\tilde{F}$ is strongly negative definite, the map $\mathcal{N}_L:\mathbb{L}^2_{\mathbb{F}}([0,T]; \mathbb{R}^d)\rightarrow \mathbb{L}^2_{\mathbb{F}}([0,T]; \mathbb{R}^d)$ is surjective and consequently:
    \begin{align*}
        \exists \lambda \in \mathbb{L}^2_{\mathbb{F}}([0,T]; \mathbb{R}^d) : \langle (I- \Phi_L)(\lambda), -\mathcal{N}_L(\lambda)\rangle_{\mathbb{L}^2_{\mathbb{F}}} = 0 
    \end{align*}
    But we have already established that it must be the case that:
    \begin{align*}
        \langle (I-\Phi_L)(\lambda), -\mathcal{N}_L(\lambda)\rangle_{\mathbb{L}^2_{\mathbb{F}}}  \geq  c^2|| \lambda||_{\mathbb{L}^2_{\mathbb{F}}}^2
    \end{align*}
    so we have that $\lambda = 0$ and since $\mathcal{N}_L$ is linear we have that $h = \mathcal{N}_L(0) = 0$. Thus by contradiction we can conclude that $(I-\Phi_L)$ is surjective and thus $\Phi_c\in \operatorname{Ran}(I-\Phi_L)$. Taking $\lambda^* \in \mathbb{L}^2_{\mathbb{F}}([0,T]; \mathbb{R}^d)$ such that $(I-\Phi_L)(\lambda^*) = \Phi_c$, we have:
    \begin{align*}
        (I-\Phi_L)(\lambda^*) = \Phi_c \iff (I-\Phi_L)(\lambda^*) - \Phi_c = 0 \iff (I-\Phi)(\lambda^*) = 0
    \end{align*}
    Thus $\lambda^* \in \operatorname{Fix}\Phi$, but since $\operatorname{Ran}\Phi\subseteq \mathcal{A}$ it must be the case that $\lambda^* \in \mathcal{A}$. The result follows from \cref{theorem:sol_iff_fixed_point}.
    \end{proof}
\end{lemma}

\subsection{A Continuation Lemma and Existence in the General Case}\label{subsec:continuation_lemma_and_existence}
We want to establish that the general coupled forward-backward system \cref{eq:fbs_general} is solvable given that the affine coupled forward-backward system of \cref{eq:fbs_affine} is solvable. To do this we can follow essentially the same strategy as when applying the method of continuation for FBSDEs (see \cite{YongPaper:1997}, \cite{MaBook:2007}). The proof rests on the proof of a suitable a priori estimate, which is then used to prove a suitable continuation lemma (for FBSDEs see \cite{MaBook:2007} Chapter 6 Lemma 2.2 and Lemma 2.3 respectively). In this section we additionally require that the global Lipschitz assumptions of \ref{ass:A5} are satisfied. First, we define the Banach space $(\mathbb{B}, ||\cdot||_{\mathbb{B}})$ such that:
\begin{align*}
    \mathbb{B}&:=L^2(\Omega, \mathcal{F}_T, \mathbb{P}; \mathbb{R}^d)\times \mathbb{L}^2_{\mathbb{F}}([0,T]; \mathbb{R}^d)\times \mathbb{L}^2_{\mathbb{F}}([0,T]; \mathbb{R}^d) \\
    ||\gamma||_{\mathbb{B}} &:= ||g_0||_{L^2_{\mathcal{F}_T}} + ||b_0||_{\mathbb{L}^2_{\mathbb{F}}} + ||f_0||_{\mathbb{L}^2_{\mathbb{F}}}, \quad \gamma \in \mathbb{B}
\end{align*}
Now consider the following family of coupled forward-backward systems indexed by $\alpha \in [0,1]$ and $\gamma \in \mathbb{B}$:
\begin{equation}\label{eq:fbs_alpha}
    \begin{aligned}
        v_t &= v_0 + \int_0^t\mathcal{N}_s^{\alpha,\gamma}(\lambda)\,ds, \quad t \in [0,T] \\
        \lambda_t &= \mathbb{E}\left[g^{\alpha,\gamma}(v_T) + \int_t^Tb^{\alpha,\gamma}(s, v_s)\,ds\ \bigg| \ \mathcal{F}_t\right], \quad \mathbb{P}\text{-a.s.}, \ t \in [0,T]
    \end{aligned}
\end{equation}
Where for $\gamma = (g_0, b_0, f_0) \in \mathbb{B}$ the coefficients of \cref{eq:fbs_alpha} are given by:
\begin{align*}
    g^{\alpha,\gamma}(v_T) &= (1-\alpha)\tilde{g}(v_T) + \alpha g(v_T) + g_0 \\
    b^{\alpha,\gamma}(t,v_t) &= (1 - \alpha)\tilde{b}(t,v_t) + \alpha b(t,v_t)  + b_0(t), \quad t \in [0,T]\\
    \mathcal{N}^{\alpha,\gamma}(\lambda) &= (1-\alpha)\tilde{\mathcal{N}}_t(\lambda) + \alpha \mathcal{N}_t(\lambda) + f_0(t), \quad t \in [0,T]
\end{align*}
where for functions $\bar{\rho}, \tilde{\rho}:[0,T]\times \Omega\times \mathbb{R}^d\rightarrow \mathbb{R}^d$ we define the corresponding Nemytskii maps $\mathcal{N}: \mathbb{L}^2_{\mathbb{F}}([0,T]; \mathbb{R}^d)\rightarrow \mathbb{U}[0,T]$ and $\tilde{\mathcal{N}}:\mathbb{L}^2_{\mathbb{F}}([0,T]; \mathbb{R}^d)\rightarrow \mathbb{L}^2_{\mathbb{F}}([0,T]; \mathbb{R}^d)$:
\begin{align*}
    \mathcal{N}_t(\lambda)(t,\omega) = \bar{\rho}(t,\omega, \lambda(t,\omega)), \quad \tilde{\mathcal{N}}_t(\lambda)(t,\omega) = \tilde{\rho}(t,\omega, \lambda(t,\omega)), \quad (t,\omega) \in [0,T]\times \Omega
\end{align*}
In addition to assuming \ref{ass:A1}, \ref{ass:A2}, \ref{ass:B1}, \ref{ass:B2} and \ref{ass:B3}, we assume that: Functions $g, \tilde{g}, b, \tilde{b}$ satisfy the global Lipschitz assumptions of \ref{ass:A5}. Consequently the results of \cref{subsec:fp_and_uniqueness} apply. Crucially, we have global Lipschitz constants $L_g, L_b, L_\rho$ and a parameter $\mu$ (for assumption \ref{ass:B3}) for the system of \cref{eq:fbs_alpha} that are independent of $(\alpha,\gamma) \in [0,1]\times \mathbb{B}$.

\bigskip\noindent Now we define the family of operators $A^{\alpha,\gamma}:\mathbb{L}^2_{\mathbb{F}}([0,T]; \mathbb{R}^d)\rightarrow \mathbb{L}^2_{\mathbb{F}}([0,T]; \mathbb{R}^d)^*$ and thus mappings $\Phi^{\alpha,\gamma}:= R^{-1}\circ A^{\alpha,\gamma}\circ \mathcal{N}^{\alpha,\gamma}$. As in previous sections, we can establish the equivalence between pairs $(v, \lambda)$ that solve \cref{eq:fbs_alpha} and $\operatorname{Fix}\Phi^{\alpha,\gamma}\subseteq \mathcal{A}$ for each $\alpha \in [0,1]$ and $\gamma \in \mathbb{B}$.
\begin{lemma} \label{lemma:fix_iff_sol_and_unique_alpha}
    Consider the family $\{\operatorname{Fix}\Phi^{\alpha, \gamma} : \alpha \in [0,1], \gamma \in \mathbb{B}\}$. Then for each $\alpha \in [0,1]$, $\gamma \in \mathbb{B}$:
    \begin{enumerate}
        \item $(v, \lambda) \in \mathbb{S}^{2,c}_{\mathbb{F}}([0,T]; \mathbb{R}^d)\times \mathbb{S}^{2}_{\mathbb{F}}([0,T]; \mathbb{R}^d)$ solves \cref{eq:fbs_alpha} iff $\lambda \in \operatorname{Fix}\Phi^{\alpha,\gamma}$
        \item If the set $\operatorname{Fix}\Phi^{\alpha,\gamma}$ is non-empty it is a singleton.
    \end{enumerate}
    \begin{proof}
        The result follows from a direct application of \cref{theorem:sol_iff_fixed_point} and \cref{theorem:fbs_unique_sol}
    \end{proof}
\end{lemma}
\begin{lemma}[A Priori Estimate] \label{lemma:key_apriori_estimate}
    Let $\alpha \in [0,1]$, $\gamma, \gamma' \in \mathbb{B}$, and let $\lambda \in \operatorname{Fix}\Phi^{\alpha, \gamma}, \ \lambda' \in \operatorname{Fix}\Phi^{\alpha, \gamma'}$. Then there exists $C > 0$ independent of $\alpha, \gamma, \gamma'$ such that: 
    \begin{align*}
        ||\lambda - \lambda'||_{\mathbb{L}^2_{\mathbb{F}}} \leq C ||\gamma -\gamma'||_{\mathbb{B}}
    \end{align*}
    \begin{proof}
        For every $\gamma = (g_0, b_0, f_0) \in \mathbb{B}$ and $0 = (0, 0, 0) \in \mathbb{B}$ we have that:
        \begin{align*}
            A^{\alpha, \gamma}(\rho)(\varphi) &= \mathbb{E}\left[(g^{\alpha, 0}(v_T) + g_0)\cdot v_T' + \int_0^T(b^{\alpha, 0}(t, v_t) + b_0(t))\cdot v_t'\,dt\right] \\
            &= \mathbb{E}\left[g^{\alpha, 0}(v_T)\cdot v_T' + \int_0^Tb^{\alpha, 0}(t, v_t)\cdot v_t'\,dt\right] + \mathbb{E}\left[g_0\cdot v_T' + \int_0^Tb_0(t)\cdot v_t'\,dt\right]\\
            &= A^{\alpha, 0}(\rho)(\varphi) + A^{\gamma}(\varphi), \quad \forall \rho, \varphi \in \mathbb{L}^2_{\mathbb{F}}([0,T]; \mathbb{R}^d)
        \end{align*}
        Thus we have the equality:
        \begin{align*}
            A^{\alpha, \gamma}(\mathcal{N}^{\alpha, \gamma}(\lambda))(\varphi) = A^{\alpha, 0}(\mathcal{N}^{\alpha,0}(\lambda) + f_0)(\varphi) + A^{\gamma}(\varphi), \quad \forall \varphi \in \mathbb{L}^2_{\mathbb{F}}([0,T]; \mathbb{R}^d), \ \forall \lambda \in \mathcal{A}
        \end{align*}
        and the duality relation:
        \begin{align*}
            A^{\alpha, \gamma}(\mathcal{N}^{\alpha, \gamma}(\lambda))(\varphi) = \langle \Phi^{\alpha, \gamma}(\lambda), \varphi\rangle_{\mathbb{L}^2_{\mathbb{F}}}, \quad \forall \varphi \in \mathbb{L}^2_{\mathbb{F}}([0,T]; \mathbb{R}^d), \quad \forall \lambda \in \mathcal{A}
        \end{align*}
        Now taking $\gamma = (g_0, b_0, f_0), \gamma' = (g_0', b_0', f_0') \in \mathbb{B}$, $\lambda \in \operatorname{Fix}\Phi^{\alpha, \gamma}$, and $\lambda' \in \operatorname{Fix}\Phi^{\alpha, \gamma'}$ we have that:
        \begin{align*}
            ||\lambda - \lambda'||_{\mathbb{L}^2_{\mathbb{F}}} &= ||\Phi^{\alpha, \gamma}(\lambda) - \Phi^{\alpha, \gamma'}(\lambda')||_{\mathbb{L}^2_{\mathbb{F}}} \\
            &= ||A^{\alpha, \gamma}(\mathcal{N}^{\alpha, \gamma}(\lambda)) - A^{\alpha, \gamma'}(\mathcal{N}^{\alpha, \gamma'}(\lambda'))||_{(\mathbb{L}^2_{\mathbb{F}})^*} \\
            &= || A^{\alpha, 0}(\mathcal{N}^{\alpha,0}(\lambda) + f_0) -  A^{\alpha, 0}(\mathcal{N}^{\alpha,0}(\lambda') + f_0') + A^{\gamma} - A^{\gamma'}||_{(\mathbb{L}^2_{\mathbb{F}})^*} \\
            &\leq || A^{\alpha, 0}(\mathcal{N}^{\alpha,0}(\lambda) + f_0) -  A^{\alpha, 0}(\mathcal{N}^{\alpha,0}(\lambda') + f_0')||_{(\mathbb{L}^2_{\mathbb{F}})^*} + || A^{\gamma} - A^{\gamma'}||_{(\mathbb{L}^2_{\mathbb{F}})^*}
        \end{align*}
        Now for $\varphi \in \{\varphi \in \mathbb{L}^2_{\mathbb{F}}([0,T]; \mathbb{R}^d) : ||\varphi||_{\mathbb{L}^2_{\mathbb{F}}}\leq 1\}$ we have for any $\rho^1, \rho^2 \in \mathbb{U}[0,T]$ that:
        \begin{align*}
            &|(A^{\alpha, 0}(\rho^1) - A^{\alpha, 0}(\rho^2))(\varphi)| \\
            &\quad\quad = \bigg|\mathbb{E}\left[(g^{\alpha, 0}(v_T^1)-g^{\alpha, 0}(v_T^2))\cdot v_T' + \int_0^T(b^{\alpha, 0}(t, v_t^1)-b^{\alpha, 0}(t, v_t^2))\cdot v_t'\,dt\right]\bigg| \\
            &\leq \mathbb{E}[||g^{\alpha, 0}(v_T^1)-g^{\alpha, 0}(v_T^2)||_2^2]^{1/2}\mathbb{E}[|| v_T'||_2^2]^{1/2} + \\
            &\quad\quad \mathbb{E}\left[\int_0^T||b^{\alpha, 0}(t, v_t^1)-b^{\alpha, 0}(t, v_t^2)||_2^2\,dt\right]^{1/2}\mathbb{E}\left[\int_0^T|| v_t'||_2^2\,dt\right]^{1/2}
        \end{align*}
        By \cref{lemma:lambda_cadlag_sm} we have:
        \begin{align*}
            \mathbb{E}[|| v_T'||_2^2]^{1/2} \leq \sqrt{T}||\varphi||_{\mathbb{L}^2_{\mathbb{F}}}, \quad \mathbb{E}\left[\int_0^T|| v_t'||_2^2\,dt\right]^{1/2} \leq T||\varphi||_{\mathbb{L}^2_{\mathbb{F}}}
        \end{align*}
        Using the global Lipschitz assumptions of \ref{ass:A5} and again applying \cref{lemma:lambda_cadlag_sm} we obtain:
        \begin{align*}
            |(A^{\alpha, 0}(\rho^1) - A^{\alpha, 0}(\rho^2))(\varphi)| &\leq \sqrt{T}L_g\mathbb{E}[||v_T^1-v_T^2||_2^2]^{1/2} + TL_b\mathbb{E}\left[\int_0^T||v_t^1- v_t^2||_2^2\,dt\right]^{1/2} \\
            &\leq (TL_g + T^2L_b)||\rho^1 - \rho^2||_{\mathbb{L}^2_{\mathbb{F}}}
        \end{align*}
        thus we obtain a Lipschitz $L_A$ constant for $A^{\alpha, 0}$ independent of $\alpha$:
        \begin{align*}
            ||A^{\alpha, 0}(\rho^1) - A^{\alpha, 0}(\rho^2)||_{(\mathbb{L}^2_{\mathbb{F}})^*} \leq L_A||\rho^1 - \rho^2||_{\mathbb{L}^2_{\mathbb{F}}}, \quad L_A = (TL_g + T^2L_b)
        \end{align*}
        Now, clearly we have that:
        \begin{align*}
             |(A^{\gamma} - A^{\gamma'})(\varphi)| &= \bigg| \mathbb{E}\left[(g_0 - g_0')\cdot v_T' + \int_0^T(b_0(t) - b_0'(t))\cdot v_t'\,dt\right]\bigg| \\
             &\leq C_T(||g_0 - g_0'||_{L^2_{\mathcal{F}_T}} + ||b_0 - b_0'||_{\mathbb{L}^2_{\mathbb{F}}}), \quad C_T := \max\{\sqrt{T}, T\}
        \end{align*}
        thus we have:
        \begin{align*}
            ||A^{\gamma} - A^{\gamma'}||_{(\mathbb{L}^2_{\mathbb{F}})^*}  &\leq C_T(||g_0 - g_0'||_{L^2_{\mathcal{F}_T}} + ||b_0 - b_0'||_{\mathbb{L}^2_{\mathbb{F}}}) \\
            &\leq C_T||\gamma - \gamma'||_{\mathbb{B}}
        \end{align*}
        Consequently, we obtain the bound:
        \begin{align*}
            ||\lambda - \lambda'||_{\mathbb{L}^2_{\mathbb{F}}} &\leq L_A(||\mathcal{N}^{\alpha, 0}(\lambda) + \mathcal{N}^{\alpha, 0}(\lambda')||_{\mathbb{L}^2_{\mathbb{F}}} + ||\gamma - \gamma||_{\mathbb{B}}) + C_T||\gamma - \gamma'||_{\mathbb{B}} \\
            &\leq C_T'\left(||\mathcal{N}^{\alpha, 0}(\lambda) + \mathcal{N}^{\alpha, 0}(\lambda')||_{\mathbb{L}^2_{\mathbb{F}}} + ||\gamma - \gamma'||_{\mathbb{B}}\right), \quad C_T' := \max(L_A, C_T)
        \end{align*}
        Now, from the duality relation we have the equality:
        \begin{align*}
            &(A^{\alpha, \gamma}(\mathcal{N}^{\alpha, \gamma}(\lambda)) - A^{\alpha, \gamma'}(\mathcal{N}^{\alpha, \gamma'}(\lambda')))(\mathcal{N}^{\alpha, \gamma}(\lambda) - \mathcal{N}^{\alpha, \gamma'}(\lambda'))  \\ 
            &= \langle \lambda - \lambda', \mathcal{N}^{\alpha, \gamma}(\lambda) - \mathcal{N}^{\alpha, \gamma'}(\lambda')\rangle_{\mathbb{L}^2_{\mathbb{F}}} 
        \end{align*}
        Expanding both sides, we have for the LHS:
        \begin{align*}
            &(A^{\alpha, \gamma}(\mathcal{N}^{\alpha, \gamma}(\lambda)) - A^{\alpha, \gamma'}(\mathcal{N}^{\alpha, \gamma'}(\lambda')))(\mathcal{N}^{\alpha, \gamma}(\lambda) - \mathcal{N}^{\alpha, \gamma'}(\lambda')) \\
            &=(A^{\alpha, 0}(\mathcal{N}^{\alpha, \gamma}(\lambda)) - A^{\alpha, 0}(\mathcal{N}^{\alpha, \gamma'}(\lambda')))(\mathcal{N}^{\alpha, \gamma}(\lambda) - \mathcal{N}^{\alpha, \gamma'}(\lambda')) \\
            &\qquad + (A^\gamma - A^{\gamma'})(\mathcal{N}^{\alpha, 0}(\lambda) - \mathcal{N}^{\alpha, 0}(\lambda')) + (A^\gamma - A^{\gamma'})(f_0 - f'_0)
        \end{align*}
        and for the RHS:
        \begin{align*}
           \langle \lambda - \lambda', \mathcal{N}^{\alpha, \gamma}(\lambda) - \mathcal{N}^{\alpha, \gamma'}(\lambda')\rangle_{\mathbb{L}^2_{\mathbb{F}}} 
            &=\langle \lambda - \lambda', \mathcal{N}^{\alpha, 0}(\lambda) - \mathcal{N}^{\alpha, 0}(\lambda')\rangle_{\mathbb{L}^2_{\mathbb{F}}} + \langle \lambda - \lambda', f_0 - f_0'\rangle_{\mathbb{L}^2_{\mathbb{F}}}
        \end{align*}
        Thus we have the equality:
        \begin{align*}
            &(A^{\alpha, 0}(\mathcal{N}^{\alpha, \gamma}(\lambda)) - A^{\alpha, 0}(\mathcal{N}^{\alpha, \gamma'}(\lambda')))(\mathcal{N}^{\alpha, \gamma}(\lambda) - \mathcal{N}^{\alpha, \gamma'}(\lambda')) \\
            &\qquad + (A^\gamma - A^{\gamma'})(\mathcal{N}^{\alpha, 0}(\lambda) - \mathcal{N}^{\alpha, 0}(\lambda')) + (A^\gamma - A^{\gamma'})(f_0 - f'_0) \\
            &= \langle \lambda - \lambda', \mathcal{N}^{\alpha, 0}(\lambda) - \mathcal{N}^{\alpha, 0}(\lambda')\rangle_{\mathbb{L}^2_{\mathbb{F}}} + \langle \lambda - \lambda', f_0 - f_0'\rangle_{\mathbb{L}^2_{\mathbb{F}}}
        \end{align*}
        rearranging we obtain:
        \begin{align*}
            &(A^{\alpha, 0}(\mathcal{N}^{\alpha, \gamma}(\lambda)) - A^{\alpha, 0}(\mathcal{N}^{\alpha, \gamma'}(\lambda')))(\mathcal{N}^{\alpha, \gamma}(\lambda) - \mathcal{N}^{\alpha, \gamma'}(\lambda')) \\
            &-\langle \lambda - \lambda', \mathcal{N}^{\alpha, 0}(\lambda) - \mathcal{N}^{\alpha, 0}(\lambda')\rangle_{\mathbb{L}^2_{\mathbb{F}}} \\
            &= \langle \lambda - \lambda', f_0 - f_0'\rangle_{\mathbb{L}^2_{\mathbb{F}}} - (A^\gamma - A^{\gamma'})(\mathcal{N}^{\alpha, 0}(\lambda) - \mathcal{N}^{\alpha, 0}(\lambda')) - (A^\gamma - A^{\gamma'})(f_0 - f'_0)
        \end{align*}
        By monotonicity of $A^{\alpha, 0}$ we have that:
        \begin{align*}
            (A^{\alpha, 0}(\mathcal{N}^{\alpha, \gamma}(\lambda)) - A^{\alpha, 0}(\mathcal{N}^{\alpha, \gamma'}(\lambda')))(\mathcal{N}^{\alpha, \gamma}(\lambda) - \mathcal{N}^{\alpha, \gamma'}(\lambda'))  \geq 0
        \end{align*}
        and we have by anti-monotonicity of $\mathcal{N}^{\alpha, 0}$ that:
        \begin{align*}
            \langle \lambda - \lambda', \mathcal{N}^{\alpha, 0}(\lambda) - \mathcal{N}^{\alpha, 0}(\lambda')\rangle_{\mathbb{L}^2_{\mathbb{F}}} \leq 0
        \end{align*}
        Consequently it must be the case that:
        \begin{align*}
            &0\leq -\langle \lambda - \lambda', \mathcal{N}^{\alpha, 0}(\lambda) - \mathcal{N}^{\alpha, 0}(\lambda')\rangle_{\mathbb{L}^2_{\mathbb{F}}}  \\
            &\leq \langle \lambda - \lambda', f_0 - f_0'\rangle_{\mathbb{L}^2_{\mathbb{F}}} - (A^\gamma - A^{\gamma'})(\mathcal{N}^{\alpha, 0}(\lambda) - \mathcal{N}^{\alpha, 0}(\lambda')) - (A^\gamma - A^{\gamma'})(f_0 - f'_0) \\
            &\leq ||\lambda - \lambda'||_{\mathbb{L}^2_{\mathbb{F}}}||f_0 - f_0'||_{\mathbb{L}^2_{\mathbb{F}}} + |(A^\gamma - A^{\gamma'})(\mathcal{N}^{\alpha, 0}(\lambda) - \mathcal{N}^{\alpha, 0}(\lambda'))| + |(A^\gamma - A^{\gamma'})(f_0 - f'_0)| \\
            &\leq ||\lambda - \lambda'||_{\mathbb{L}^2_{\mathbb{F}}}||\gamma - \gamma||_{\mathbb{B}} + ||(A^\gamma - A^{\gamma'}||_{(\mathbb{L}^2_{\mathbb{F}})^*}(||\mathcal{N}^{\alpha, 0}(\lambda) - \mathcal{N}^{\alpha, 0}(\lambda')||_{\mathbb{L}^2_{\mathbb{F}}} + ||\gamma - \gamma'||_{\mathbb{B}})
        \end{align*}
        From our assumption \ref{ass:B3}, and recalling that $||A^{\gamma} - A^{\gamma'}||_{(\mathbb{L}^2_{\mathbb{F}})^*} \leq C_T||\gamma - \gamma'||_{\mathbb{B}}$, we obtain:
        \begin{align*}
            &\mu||\mathcal{N}^{\alpha, 0}(\lambda) - \mathcal{N}^{\alpha, 0}(\lambda')||_{_{\mathbb{L}^2_{\mathbb{F}}}}^2 \leq ||\lambda - \lambda'||_{\mathbb{L}^2_{\mathbb{F}}}||\gamma - \gamma||_{\mathbb{B}} \\
            &\quad\quad + C_T||\gamma - \gamma'||_{\mathbb{B}}||\mathcal{N}^{\alpha, 0}(\lambda) - \mathcal{N}^{\alpha, 0}(\lambda')||_{\mathbb{L}^2_{\mathbb{F}}} + C_T||\gamma - \gamma'||_{\mathbb{B}}^2
        \end{align*}
        Now recalling the earlier bound we derived for $||\lambda - \lambda'||_{\mathbb{L}^2_{\mathbb{F}}}$:
        \begin{align*}
            ||\lambda - \lambda'||_{\mathbb{L}^2_{\mathbb{F}}} \leq C_T'\left(||\mathcal{N}^{\alpha, 0}(\lambda) + \mathcal{N}^{\alpha, 0}(\lambda')||_{\mathbb{L}^2_{\mathbb{F}}} + ||\gamma - \gamma'||_{\mathbb{B}}\right)
        \end{align*}
        We have:
        \begin{align*}
            &\mu||\mathcal{N}^{\alpha, 0}(\lambda) - \mathcal{N}^{\alpha, 0}(\lambda')||_{_{\mathbb{L}^2_{\mathbb{F}}}}^2 \leq 2C_T'\left(||\mathcal{N}^{\alpha, 0}(\lambda) + \mathcal{N}^{\alpha, 0}(\lambda')||_{\mathbb{L}^2_{\mathbb{F}}}||\gamma - \gamma||_{\mathbb{B}} + ||\gamma - \gamma'||_{\mathbb{B}}^2\right)
        \end{align*}
        Now by Young's inequality we have for $\epsilon > 0$ that:
        \begin{align*}
            ||\mathcal{N}^{\alpha, 0}(\lambda) + \mathcal{N}^{\alpha, 0}(\lambda')||_{\mathbb{L}^2_{\mathbb{F}}}||\gamma - \gamma||_{\mathbb{B}}  \leq \frac{\epsilon}{2}||\mathcal{N}^{\alpha, 0}(\lambda) + \mathcal{N}^{\alpha, 0}(\lambda')||_{\mathbb{L}^2_{\mathbb{F}}}^2 + \frac{1}{2\epsilon}||\gamma - \gamma||_{\mathbb{B}}^2 
        \end{align*}
        Choosing $\epsilon = \mu/(2C_T')$ we obtain:
        \begin{align*}
            &||\mathcal{N}^{\alpha, 0}(\lambda) - \mathcal{N}^{\alpha, 0}(\lambda')||_{{\mathbb{L}^2_{\mathbb{F}}}}^2 \leq C_{T,\mu}||\gamma - \gamma||_{\mathbb{B}}^2, \quad C_{T,\mu} := \frac{2C_T'}{\mu}\max\left(\frac{1}{\mu}, 2\right)
        \end{align*}
        Now feeding this into our existing bound for $||\lambda - \lambda'||_{\mathbb{L}^2_{\mathbb{F}}}$, we obtain the result:
        \begin{align*}
            ||\lambda - \lambda'||_{\mathbb{L}^2_{\mathbb{F}}} \leq C||\gamma - \gamma||_{\mathbb{B}}
        \end{align*}
    \end{proof}
\end{lemma}
\begin{lemma}[Continuation Lemma]\label{lemma:continuation_lemma}
    Suppose for some $\alpha \in [0,1)$ and $\forall \gamma \in \mathbb{B}$ we have that $\operatorname{Fix}\Phi^{\alpha,\gamma}\neq\emptyset$, then there exists $\delta > 0$ with $\alpha + \delta \leq 1$ such that $\operatorname{Fix}\Phi^{\alpha+\epsilon,\gamma}\neq\emptyset$ for every $\epsilon \in (0,\delta]$.
    \begin{proof}
        We use essentially the same strategy as Lemma 2.3 in Chapter 6 of \cite{MaBook:2007}. With $\delta > 0$ to be determined, let $0< \epsilon \leq \delta$ be given. Given a sequence $(\gamma_k)_{k \in \mathbb{N}}\subseteq \mathbb{B}$ we have by assumption that $\operatorname{Fix}\Phi^{\alpha, \gamma_k}\neq \emptyset$. Now taking $\gamma_0 = (g_0, b_0, f_0) \in \mathbb{B}$, for $k \geq 1$ we can define the sequence of operators $\Phi^{\alpha, \gamma_k}$. Now choosing $\lambda^0, v^0 = 0$, for $k\geq 1$ we choose $\lambda^{k} \in \operatorname{Fix}\Phi^{\alpha, \gamma_{k}}$ where by \cref{lemma:fix_iff_sol_and_unique_alpha} we can identify  $\lambda^k \in \operatorname{Fix}\Phi^{\alpha,\gamma_k}$ with a solution $(v^k, \lambda^k)$ to the equation:
        \begin{align*}
            v_t^k &= v_0 + \int_0^t\mathcal{N}^{\alpha, \gamma_k}_s(\lambda^k)\,ds \quad t \in [0,T] \\
            \lambda^k_t &= \mathbb{E}\left[g^{\alpha, \gamma_k}(v^k_T) + \int_t^Tb^{\alpha, \gamma_k}(t, v^k_t)\,ds \ \bigg| \mathcal{F}_t\right], \quad \mathbb{P}\text{-a.s.}, \quad t \in [0,T]
        \end{align*}
        which is simply \cref{eq:fbs_alpha} with particular coefficients defined by:
        \begin{align*}
            g^{\alpha, \gamma_k}(v) &:= (1-\alpha)\tilde{g}(v) + \alpha g(v) + g_0^{k} \\
            b^{\alpha,\gamma_k}(t,v) &:= (1 - \alpha)\tilde{b}(t,v) + \alpha b(t,v) + b_0^{k}(t) \\
            \mathcal{N}^{\alpha,\gamma_k}(\lambda) &:= (1-\alpha)\tilde{\mathcal{N}}(\lambda) + \alpha \mathcal{N}(\lambda) + f_0^{k}
        \end{align*}
        where we choose the sequence $(\gamma_k = (g_0^k, b_0^k, f_0^k))_{k \in \mathbb{N}}\subseteq \mathbb{B}$ as follows:
        \begin{align*}
            g_0^{k} &= \epsilon(g(v_T^{k-1}) - \tilde{g}(v_T^{k-1})) + g_0 \\
            b_0^{k}(t) &= \epsilon(b(t,v_t^{k-1}) - \tilde{b}(t,v_t^{k-1})) + b_0(t), \quad  t \in [0,T]\\
            f_0^{k}(t) &= \epsilon(\mathcal{N}_t(\lambda^{k-1}) - \tilde{\mathcal{N}}_t(\lambda^{k-1})) + f_0(t), \quad t \in [0,T]
        \end{align*}
         By \cref{lemma:key_apriori_estimate} we have that there exists $C_1 > 0$ such that:
        \begin{align*}
            ||\lambda^{k+1} - \lambda^{k}||_{\mathbb{L}^2_{\mathbb{F}}} \leq C_1||\gamma^{k+1} - \gamma^k||_{\mathbb{B}}, \quad k = 1, 2, \dots
        \end{align*}
        Recalling that:
        \begin{align*}
            &||\gamma_{k+1} - \gamma_k||_{\mathbb{B}} := ||g_0^{k+1} -g_0^{k}||_{L^2_{\mathcal{F}_T}} +  ||b_0^{k+1}- b_0^{k}||_{\mathbb{L}^2_{\mathbb{F}}} + ||f_0^{k+1}- f_0^{k}||_{\mathbb{L}^2_{\mathbb{F}}}, \quad k= 1, 2, \dots 
        \end{align*}
        by the global Lipschitz conditions of \ref{ass:A5}, for $k\geq 1$ we have that:
        \begin{align*}
            ||g_0^{k+1} -g_0^{k}||_{L^2_{\mathcal{F}_T}} &\leq \epsilon\left(||g(v_T^{k}) -g(v_T^{k-1})||_{L^2_{\mathcal{F}_T}} + ||\tilde{g}(v_T^{k})) -\tilde{g}(v_T^{k-1})) ||_{L^2_{\mathcal{F}_T}}\right) \\ 
            &\leq \epsilon 2L_g||v_T^{k} -v_T^{k-1}||_{L^2_{\mathcal{F}_T}} \\ \\
            ||b_0^{k+1}- b_0^{k}||_{\mathbb{L}^2_{\mathbb{F}}} &\leq \epsilon\mathbb{E}\left[\int_0^T||b(t,v_t^k) - b(t,v_t^{k-1})||_2^2\,dt\right]^{1/2} \\ &\quad\quad\quad+ \epsilon \mathbb{E}\left[\int_0^T||\tilde{b}(t,v_t^k) - \tilde{b}(t,v_t^{k-1})||_2^2\,dt\right]^{1/2} \\
            &\leq \epsilon 2L_b||v^k - v^{k-1}||_{\mathbb{L}^2_{\mathbb{F}}} \\\\
            ||f_0^{k+1}- f_0^{k}||_{\mathbb{L}^2_{\mathbb{F}}} &\leq \epsilon \left(||\mathcal{N}(\lambda^k) - \mathcal{N}(\lambda^{k-1})||_{\mathbb{L}^2_{\mathbb{F}}} + ||\tilde{\mathcal{N}}(\lambda^k) - \tilde{\mathcal{N}}(\lambda^{k-1})||_{\mathbb{L}^2_{\mathbb{F}}}\right) \\
            &\leq \epsilon 2L_\rho||\lambda^k - \lambda^{k-1}||_{\mathbb{L}^2_{\mathbb{F}}}
        \end{align*}
        So we can conclude that there exists $C_2 > 0$ such that:
        \begin{align*}
            ||\gamma_{k+1} - \gamma_{k}||_{\mathbb{L}^2_{\mathbb{F}}} \leq \epsilon C_2\left(||\lambda^k - \lambda^{k-1}||_{\mathbb{L}^2_{\mathbb{F}}} + ||v^{k} - v^{k-1}||_{\mathbb{S}^{2,c}_{\mathbb{F}}}\right), \quad k\geq 1
        \end{align*}
        Consequently we have that:
        \begin{align*}
            ||\lambda^{k+1} - \lambda^{k}||_{\mathbb{L}^2_{\mathbb{F}}} \leq \epsilon C_1C_2\left(||\lambda^k - \lambda^{k-1}||_{\mathbb{L}^2_{\mathbb{F}}} + ||v^{k} - v^{k-1}||_{\mathbb{S}^{2,c}_{\mathbb{F}}}\right), \quad k\geq 1
        \end{align*}
        Now, we have for $t \in [0,T]$ that:
        \begin{align*}
            ||v^{k+1} - v^{k}||_{\mathbb{S}^{2,c}_{\mathbb{F}}} &\leq \sqrt{T}||\mathcal{N}^{\alpha, \gamma_{k+1}}(\lambda^{k+1}) - \mathcal{N}^{\alpha, \gamma_{k}}(\lambda^{k})||_{\mathbb{L}^2_{\mathbb{F}}} \\
            &= \sqrt{T}||\mathcal{N}^{\alpha, 0}(\lambda^{k+1}) - \mathcal{N}^{\alpha, 0}(\lambda^{k}) + f_0^{k+1} - f_0^k||_{\mathbb{L}^2_{\mathbb{F}}} \\
            &\leq \sqrt{T}L_\rho||\lambda^{k+1} - \lambda^k||_{\mathbb{L}^2_{\mathbb{F}}} + \sqrt{T}||f_0^{k+1} - f_0^k||_{\mathbb{L}^2_{\mathbb{F}}} \\
            &\leq \sqrt{T}L_\rho||\lambda^{k+1} - \lambda^k||_{\mathbb{L}^2_{\mathbb{F}}} + \sqrt{T}||\gamma_{k+1} - \gamma_k||_{\mathbb{B}}
        \end{align*}
        Thus we have that:
        \begin{align*}
             ||v^{k+1} - v^{k}||_{\mathbb{S}^{2,c}_{\mathbb{F}}} \leq \epsilon\sqrt{T}(L_pC_1C_2 + C_2)\left(||\lambda^k - \lambda^{k-1}||_{\mathbb{L}^2_{\mathbb{F}}} + ||v^{k} - v^{k-1}||_{\mathbb{S}^{2,c}_{\mathbb{F}}}\right) 
        \end{align*}
        Thus we finally obtain the inequality:
        \begin{align*}
             ||v^{k+1} - v^{k}||_{\mathbb{S}^{2,c}_{\mathbb{F}}} + ||\lambda^{k+1} - \lambda^{k}||_{\mathbb{L}^2_{\mathbb{F}}} \leq \epsilon C\left(||v^{k} - v^{k-1}||_{\mathbb{S}^{2,c}_{\mathbb{F}}} + ||\lambda^k - \lambda^{k-1}||_{\mathbb{L}^2_{\mathbb{F}}}\right), \quad k\geq 1
        \end{align*}
        Where we have: $C := (C_1C_2(1 + \sqrt{T}L_p) + \sqrt{T}C_2)$. Now choosing $\delta$ such that $0 < \delta < 1/C$ and $\delta\leq 1-\alpha$, we have that the sequence $(v^k, \lambda^k)_{k\geq 1}$ is Cauchy in the Banach space $\mathbb{S}^{2,c}_{\mathbb{F}}([0,T]; \mathbb{R}^d)\times \mathbb{L}^2_{\mathbb{F}}([0,T]; \mathbb{R}^d)$ and thus converges to a unique limit point $(v^*, \lambda^*)$. Recall that by our construction of the iteration procedure and the sequence $(\gamma_k)_{k \in \mathbb{N}}\subseteq \mathbb{B}$, it must be the case that:
        \begin{align*}
            v^*_t &= v_0 + \int_0^t\mathcal{N}^{\alpha + \epsilon, \gamma_0}_s(\lambda^*)\,ds, \quad t \in [0,T] \\
            \lambda^*_t &= \mathbb{E}\left[g^{\alpha+\epsilon, \gamma_0}(v^*_T) + \int_t^Tb^{\alpha+\epsilon, \gamma_0}(s, v_s^*)\,ds \ \bigg| \mathcal{F}_t\right], \quad \mathbb{P}\text{-a.s.}, \quad t \in [0,T]
        \end{align*}
        Consequently $\lambda^* \in \operatorname{Fix}\Phi^{\alpha+\epsilon,\gamma_0}$ by \cref{lemma:fix_iff_sol_and_unique_alpha}, and since $\gamma_0 \in \mathbb{B}$ was arbitrary we have that $\operatorname{Fix}\Phi^{\alpha+\epsilon,\gamma}\neq \emptyset$ for each $\gamma \in \mathbb{B}$ for every $\epsilon \in (0,\delta]$.
    \end{proof}
\end{lemma}
\begin{theorem}[Existence of Fixed Points] \label{theorem:fbs_general_existence}
    Defining $\Phi$ as in \cref{subsec:fp_and_uniqueness}. The set $\operatorname{Fix}\Phi\neq\emptyset$, and consequently a solution to the coupled forward backward system of \cref{eq:fbs_general} exists.
    \begin{proof}
        Let $\tilde{G}, \tilde{B}, \tilde{F} \in \mathbb{R}^{d\times d}$ satisfy the assumptions of \cref{subsec:existence_affine_system}, then with 
        \begin{align*}
            \tilde{g}(\omega, v) = \tilde{G}v, \quad \tilde{b}(\omega, t,v) = \tilde{B}v, \quad \tilde{\rho}(t,\omega,\lambda) = \tilde{F}\lambda, \quad (t,\omega, v, \lambda) \in [0,T]\times \Omega\times \mathbb{R}^d\times \mathbb{R}^d
        \end{align*}
        then $\forall \gamma = (g_0, b_0, f_0) \in \mathbb{B}$, the affine system:
        \begin{align*}
            v_t &= v_0 + \int_0^t(\tilde{\rho}(s,\lambda_s) + f_0(s))\,ds, \quad t \in [0,T] \\
            \lambda_t &= \mathbb{E}\left[\tilde{g}(v_T) + \int_0^T(\tilde{b}(s,v_s) + b_0(s))\,ds \ \bigg|\mathcal{F}_t\right], \quad \mathbb{P}\text{-a.s.}, \quad t \in [0,T]
        \end{align*}
        has a unique solution by \cref{lemma:uniqueness_affine_case} and \cref{lemma:existence_affine_case}. We can now define the family of coupled forward backward systems as in \cref{eq:fbs_alpha} parameterized by $(\alpha,\gamma) \in [0,1]\times \mathbb{B}$ related to the family of operators $\{\Phi^{\alpha,\gamma} : \alpha \in [0,1], \gamma \in \mathbb{B}\}$. Our affine system corresponds to the family $\{\Phi^{0,\gamma}; \gamma \in \mathbb{B}\}$ where we know that $\operatorname{Fix}\Phi^{0,\gamma}\neq\emptyset,\ \forall \gamma \in \mathbb{B}$, thus by a finite number of applications of \cref{lemma:continuation_lemma} we have that $\operatorname{Fix}\Phi^{1,\gamma}\neq\emptyset,\ \forall \gamma \in \mathbb{B}$. Since the operator $\Phi^{1,\gamma}$ corresponds to the system:
         \begin{align*}
            v_t &= v_0 + \int_0^t(\bar{\rho}(s,\lambda_s) + f_0(s))\,ds, \quad t \in [0,T] \\
            \lambda_t &= \mathbb{E}\left[g(v_T) + \int_t^T(b(s,v_s) + b_0(s))\,ds \ \bigg|\mathcal{F}_t\right], \quad \mathbb{P}\text{-a.s.}, \quad t \in [0,T]
        \end{align*}
        with $\gamma = (g_0, b_0, f_0) \in \mathbb{B}$. Taking $\lambda^* \in \operatorname{Fix}\Phi^{1,\gamma}$, by \cref{lemma:fix_iff_sol_and_unique_alpha} the pair $(v^*,\lambda^*) \in \mathbb{S}^{2,c}_{\mathbb{F}}([0,T]; \mathbb{R}^d)\times \mathbb{S}^2_{\mathbb{F}}([0,T]; \mathbb{R}^d)$ is a solution. The result follows by taking $\gamma = (0,0,0)$ and noting that $\Phi = \Phi^{1,0}$.
    \end{proof}
\end{theorem}

\section{Conclusion}

In this paper we have taken a functional analytic approach to proving existence and uniqueness of solutions for a certain class of coupled forward-backward systems defined on a general probability space assumed to satisfy only the usual conditions. By relating a large set of decoupled forward-backward systems (containing the coupled system) to a certain strictly monotone operator via an important duality pairing, we discover we can associate solutions to the coupled equation to solutions to a fixed point problem on a suitable Hilbert space. Thus questions of existence and uniqueness of solutions to the coupled forward-backward system can be answered by studying the unique solvability of the fixed point problem. This conceptual re-framing turns out to be very fruitful when it comes to proving the required but fairly difficult a-priori estimate and continuation lemma from which we obtain an existence proof.

\bibliographystyle{alpha}
\bibliography{references}

\appendix

\section{Auxiliary Results}
\begin{definition}
    Define the set of simple $\mathbb{F}$-adapted essentially bounded processes that are c\`adl\`ag (resp. c\`agl\`ad):
\begin{align*}
    \mathcal{S}^{b, c\`adl\`ag}_{\mathbb{F}}([0,T]; \mathbb{R}^d):=\{\varphi:[0,T]\times \Omega \mapsto \mathbb{R}^d: \varphi(t,\omega) := \sum_{i=1}^{n-1}\xi_i(\omega)1_{[t_i, t_{i+1})}(t) + \xi_N(\omega)1_{\{T\}}(t) , \\  \xi_i \in L^\infty(\Omega, \mathcal{F}_{t_i}, \mathbb{P}; \mathbb{R}^d), \ 0 = t_1 < t_2 < \dots, t_n = T, \ n \in \mathbb{N}\}
\end{align*}
\begin{align*}
    \mathcal{S}^{b, c\`agl\`ad}_{\mathbb{F}}([0,T]; \mathbb{R}^d):=\{\varphi:[0,T]\times \Omega \mapsto \mathbb{R}^d: \varphi(t,\omega) := \xi_0(\omega)1_{\{0\}}(t) + \sum_{i=1}^{n-1}\xi_i(\omega)1_{(t_i, t_{i+1}]}(t), \\  \xi_i \in L^\infty(\Omega, \mathcal{F}_{t_i}, \mathbb{P}; \mathbb{R}^d), \ 0 = t_1 < t_2 < \dots, t_n = T, \ n \in \mathbb{N}\}
\end{align*}    
\end{definition}
\begin{lemma}\label{lemma:simple_proc_dense_LpF}
    The sets $\mathcal{S}^{b, c\`adl\`ag}_{\mathbb{F}}([0,T]; \mathbb{R}^d)$ and $\mathcal{S}^{b, c\`agl\`ad}_{\mathbb{F}}([0,T]; \mathbb{R}^d)$ are dense in $\mathbb{L}^p_{\mathbb{F}}([0,T]; \mathbb{R}^d)$ for $p \in \mathbb{N}$.
    \begin{proof}
        Recall that the set $\mathcal{S}^b_{\mathbb{F}}$ is dense in $\mathbb{L}^p_{\mathbb{F}}([0,T]; \mathbb{R}^d)$ iff for every $X \in \mathbb{L}^p_{\mathbb{F}}([0,T]; \mathbb{R}^d)$ there exists $\{X_n\}_{n \in \mathbb{N}}\subseteq \mathcal{S}^b_{\mathbb{F}}$ with $X_n\rightarrow X$. Such a sequence is constructed in the case of $\mathbb{L}^2_{\mathbb{F}}([0,T]; \mathbb{R})$ in \cite{ShreveBook:1996} 3.A lemma 2.4 and proposition 2.6. The construction can be easily modified to construct an approximating sequence for any $X \in \mathbb{L}^p_{\mathbb{F}}([0,T]; \mathbb{R}^d)$.
    \end{proof}
\end{lemma}
\begin{lemma}\label{lemma:firm_non_expans}
    Let $A \in \mathbb{R}^{d\times d}$ be strongly positive definite, then there exists $\mu > 0$ such that:
    \begin{align*}
        \langle Ax, x\rangle \geq \mu||Ax||_2^2, \quad \forall x \in \mathbb{R}^d
    \end{align*}
    \begin{proof}
        We have that:
        \begin{align*}
            ||Ax||_2^2 \leq ||A||_{op}^2||x||_2^2, \quad \forall x \in \mathbb{R}^d
        \end{align*}
        where $||A||_{op} > 0$, and since $A$ is strongly positive definite we have for some $c > 0$:
        \begin{align*}
            \langle Ax, x\rangle \geq c||x||_2^2, \quad \forall x \in \mathbb{R}^d
        \end{align*}
        Thus we have:
        \begin{align*}
            \langle Ax, x\rangle \geq c||x||_2^2 \geq \frac{c}{||A||_{op}}||Ax||_2^2, \quad \forall x \in \mathbb{R}^d
        \end{align*}
        so the result follows with $\mu = c/||A||_{op}$.
    \end{proof}
\end{lemma}
\end{document}